\documentstyle{amsppt}
\baselineskip18pt
\magnification=\magstep1
\pagewidth{30pc}
\pageheight{45pc}
\hyphenation{co-deter-min-ant co-deter-min-ants pa-ra-met-rised
pre-print pro-pa-gat-ing pro-pa-gate
fel-low-ship Cox-et-er dis-trib-ut-ive}
\def\leaderfill{\leaders\hbox to 1em{\hss.\hss}\hfill}
\def\A{{\Cal A}}
\def\D{{\Cal D}}
\def\H{{\Cal H}}

\def\Sy{{\Cal S}}

\def\afn{{\text {\bf a}}}

\def\idest{i.e.,\ }
\def\wh{\widehat}

\def\g{{\gamma}}
\def\G{{\Gamma}}
\def\d{{\delta}}

\def\e{{\varepsilon}}

\def\k{{\kappa}}
\def\l{{\lambda}}

\def\t{{\tau}}

\def\ba{{\bold a}}

\def\bc{{\bold c}}

\def\BB{{\bold B}}
\def\b0{\text{\bf 0}}

\def\ra{{\ \longrightarrow \ }}

\def\supp{\text{\rm \, supp}}

\def\lan{{\langle}}
\def\ran{{\rangle}}

\def\extln{{{\Cal D}(\widehat A_{n-1})}}
\def\extll{{{\Cal D}(\widehat A_l)}}
\def\tln{{TL(\widehat A_n)}}

\def\annn{\text{\rm Ann({\bf n})}}

\def\dt{{\Bbb D \Bbb T}}
\def\real{{\Bbb R}}
\def\complex{{\Bbb C}}
\def\zed{{\Bbb Z}}
\def\kyu{{\Bbb Q}}
\def\enn{{\Bbb N}}

\def\Im{\text{\rm Im}}

\def\boxit#1{\vbox{\hrule\hbox{\vrule \kern3pt
\vbox{\kern3pt\hbox{#1}\kern3pt}\kern3pt\vrule}\hrule}}
\def\rabbit{\vbox{\hbox{\kern0pt
\vbox{\kern0pt{\hbox{---}}\kern3.5pt}}}}

\def\tableau#1{
        \hbox {
                \hskip -10pt plus0pt minus0pt
                \raise\baselineskip\hbox{
                \offinterlineskip
                \hbox{#1}}
                \hskip0.25em
        }
}

\def\tabCol#1{
\hbox{\vtop{\hrule
\halign{\strut\vrule\hskip0.5em##\hskip0.5em\hfill\vrule\cr\lower0pt
\hbox\bgroup$#1$\egroup \cr}
\hrule
} } \hskip -10.5pt plus0pt minus0pt}

\def\CR{
        $\egroup\cr
        \noalign{\hrule}
        \lower0pt\hbox\bgroup$
}



\def\blank#1#2{
\hbox to #1{\hfill \vbox to #2{\vfill}}
}


\def\strut{\vrule height10pt depth5pt width0pt}

\topmatter
\title Tabular algebras and their asymptotic versions
\endtitle

\author R.M. Green \endauthor
\affil 
Department of Mathematics and Statistics\\ Lancaster University\\
Lancaster LA1 4YF\\ England\\
{\it  E-mail:} r.m.green\@lancaster.ac.uk\\
\endaffil

\abstract
We introduce tabular algebras, which are simultaneous generalizations
of cellular algebras (in the sense of Graham--Lehrer) and table 
algebras (in the sense of Arad--Blau).  We show that if a tabular
algebra is equipped with a certain kind of trace map then the algebra
has a corresponding asymptotic version whose structure can be
explicitly determined.  We also study various natural examples of
tabular algebras.
\endabstract

\thanks
The author was supported in part by an award from the Nuffield
Foundation.
\endthanks


\endtopmatter

\centerline{\bf To appear in the Journal of Algebra}

\head Introduction \endhead

The purpose of this paper is to introduce tabular algebras, a
class of associative algebras over $\zed[v, v^{-1}]$.  A tabular
algebra is defined with a ``tabular'' basis and an
anti-automorphism which are required to satisfy various properties.
The construction is a simultaneous generalization of table algebras and of
cellular algebras.  Our primary objects of study are the tabular
bases, and we find that there are important examples of associative
algebras in the literature equipped with ``natural'' bases that turn
out to be tabular.

Table algebras were introduced by Arad and Blau \cite{{\bf 1}} in order to
study irreducible characters and conjugacy classes of finite groups in
an abstract setting.  Table algebras are related in a precise way to
the association schemes of algebraic combinatorics \cite{{\bf 3}} 
and to Kawada's $C$-algebras \cite{{\bf 17}}.  The table algebras in this
paper are more general and similar to the discrete hypergroups
appearing in the work of Sunder \cite{{\bf 29}}; these have applications to
subfactors.

Cellular algebras were introduced by Graham and Lehrer \cite{{\bf 8}},
and are a class of finite dimensional associative algebras defined in
terms of a ``cell datum'' and three axioms.  One of their main
strengths is that it is relatively straightforward to construct and to
classify the irreducible modules for a cellular algebra.  Theorem
2.1.1 gives a useful and sufficient criterion for a tabular algebra 
to be cellular.

There is another---completely different---definition of
``cellular algebra'' in the literature which is 
due to Lehman and Weisfeiler \cite{{\bf 21}}
and which, ironically, is closely related to association schemes.  We
are not concerned with these algebras here.

We also introduce the notion of a ``tabular algebra with trace''.  In
this situation, the tabular algebra is equipped with an $\ba$-function 
analogous to the $\ba$-function appearing in Lusztig's work 
\cite{{\bf 22}, {\bf 23}, {\bf 24}, {\bf 25}} 
and a trace map which is compatible with this
$\ba$-function in a certain sense which will be made precise.  Tabular
algebras with trace have some interesting properties which we
investigate, such as the existence of a bilinear form (Theorem 2.2.5) 
that makes transparent the structure of the algebra as a symmetric algebra.
We also define asymptotic analogues of tabular algebras
by sending the parameter $v$ to $\infty$ in a controlled way
using a general method due to Lusztig \cite{{\bf 26}} which makes the
structure of the algebra over a suitable field explicit (Theorem 3.2.4).

The second half of the paper is devoted to the detailed study of
certain classes of examples of tabular
algebras.  Our motivation for studying these objects comes from the
canonical bases for generalized Temperley--Lieb algebras introduced by
the author and J. Losonczy \cite{{\bf 15}}.  The latter often give, or are
closely related to, examples of tabular algebras with trace.
There are other interesting examples which we mention in less
detail, such as the Hecke algebra of type $A$, Jones' annular
algebra and the Brauer algebra.  We conclude with some questions.

\head 1. Preliminaries \endhead

In \S1, we recall the definitions of table algebras and cellular
algebras, and show how they may be used to introduce the notion of a
tabular algebra.  Table algebras will always be defined over subrings of
$\complex$, typically $\zed$, and tabular algebras will always be
defined over the ring of Laurent polynomials, $\zed[v, v^{-1}]$.  For
various purposes it is convenient to extend scalars by suitable
tensoring; this will usually be made explicit.

\subhead 1.1 Table algebras \endsubhead

We begin by defining table algebras, which were introduced in the
finite dimensional commutative case by Arad and Blau \cite{{\bf 1}}.

\definition{Definition 1.1.1}
A table algebra is a pair $(A, \BB)$, where $A$ is an associative 
unital $R$-algebra for some $\zed \leq R \leq \complex$ 
and $\BB = \{b_i : i \in I\}$ is a distinguished basis for
$A$ such that $1 \in \BB$, satisfying the following three axioms:

\item{(T1)}{The structure constants of $A$ with respect to the basis
$\BB$ lie in $\real^+$, the nonnegative real numbers.}
\item{(T2)}{There is an algebra anti-automorphism $\bar{\ }$ of $A$ whose
square is the identity and which has the property 
that $b_i \in \BB \Rightarrow \overline{b_i} \in
\BB$.  (We define $\overline{i}$ by the condition $\overline{b_i} =
b_{\bar{i}}$.)}
\item{(T3)}{Let $\k(b_i, a)$ be the coefficient of $b_i$ in $a \in A$.
Then there is a function $g: \BB \times \BB \ra \real^+$ satisfying $$
\k(b_m, b_i b_j) = g(b_i, b_m) \k(b_i, b_m \overline{b_j})
,$$ where $g(b_i, b_m)$ is independent of $j$, for all $i, j, m$.
}
\enddefinition

If $R = \complex$, it follows from \cite{{\bf 1}, Proposition 2.2} 
that the basis elements of
a finite dimensional, commutative table algebra may be uniquely
rescaled so that the function $g$ in
axiom (T3) sends all pairs of basis elements to $1 \in \real$.  This
motivates the following

\definition{Definition 1.1.2}
A {\it normalized} table algebra $(A, \BB)$ over $R$ is one whose structure
constants lie in $\zed$ and for which the function $g$ in axiom (T3)
sends all pairs of basis elements to $1 \in \real$.  All table
algebras from now on will be normalized.
\enddefinition

\remark{Remark 1.1.3}
Our definition is more general than Arad--Blau's original
definition in the sense that we allow $A$ to be noncommutative and/or infinite 
dimensional.  There are many variants of the definition in the
literature, such as the generalized table algebras of Arad, Fisman and
Muzychuk \cite{{\bf 2}}, which are not required to be commutative but are
of finite rank.  A normalized table algebra (in our sense) 
corresponds to a discrete hypergroup in the sense of
Sunder \cite{{\bf 29}, Definition IV.1}.  The table algebras in this
paper are typically finite dimensional and typically commutative.
\endremark

\definition{Definition 1.1.4}
If $(A, \BB)$ is a table algebra and $a \in A$, we write $\supp(a)$
to denote the set of elements of $\BB$ which occur with nonzero
coefficient in $A$.  (Arad and Blau use the notation $\text{Irr}(a)$.)
\enddefinition

A table algebra may be viewed geometrically, as follows.

\proclaim{Lemma 1.1.5}
Let $(A, \BB)$ be a normalized table algebra over $R \leq \real$ with
linear anti-automorphism $\bar{\ }$.  
There exists a unique positive definite symmetric bilinear form 
$h : (\real \otimes_R A) \times (\real \otimes_R A) \ra \real$ 
with the property that $h(ab, c) = h(a, c\overline{b})$
for all $a, b, c \in A$ and with respect to which the set $\BB$ is an 
orthonormal basis.
\endproclaim

\demo{Proof}
If $b, b' \in \BB$, the orthonormal basis hypothesis requires 
$h(b, b') = \delta_{b, b'}$, so we
define the function $h$ to be the unique (symmetric) bilinear form with
this property.  The fact that $h$ is positive definite follows
easily.  To check
$h(ab, c) = h(a, c\overline{b})$ for all $a, b, c \in
A$, it is enough to consider basis elements $a, b, c$.
However, the definition of $f$, axioms (T1), (T2), (T3) and the fact that
the algebra is normalized imply that $$
h(ab, c) = \kappa(c, ab) = \kappa(a, c\overline{b}) = h(a, c\overline{b})
,$$ as required.
\qed\enddemo

The next result was proved by Arad and Blau \cite{{\bf 1}, Proposition 2.5}
in the commutative case, and is related to \cite{{\bf 2}, Corollary 3.5}.

\proclaim{Proposition 1.1.6}
Let $(A, \BB)$ be a finite dimensional normalized table algebra over
$R$ (where $\zed \leq R \leq \real$) with linear anti-automorphism
$\bar{\ }$.  Then $\real \otimes_R A$ is semisimple as an $\real$-algebra.
\endproclaim

\demo{Proof}
It suffices to prove that $\real \otimes_R A$ is semisimple as a right 
module over itself, or that any right ideal in $A$ is complemented.  
Let $K$ be a right ideal of $\real \otimes_R A$, and let $K'$ be the 
orthogonal complement of $K$ in $\real \otimes_R A$ with respect to 
the inner product $h$ of Lemma 1.1.5.
Consideration of the equation $h(ab, c) = h(a,
c\overline{b})$ where $a \in K$, $c \in K'$ shows that $K'$ is also a
right ideal, which completes the proof.
\qed\enddemo
\subhead 1.2 Cellular algebras \endsubhead

Cellular algebras were originally defined by Graham and Lehrer
\cite{{\bf 8}}.  Like table algebras, they are associative algebras with
involution which are defined to satisfy certain axioms.

\definition{Definition 1.2.1}
Let $R$ be a commutative
ring with identity.  A {\it cellular algebra} over $R$ is an associative unital
algebra, $A$, together with a cell datum $(\Lambda, M, C, *)$ where:

\item {(C1)}
{$\Lambda$ is a finite poset.  For each $\l \in \Lambda$, $M(\l)$ is a 
finite set (the set of ``tableaux'' of type $\l$) such that $$
C : \coprod_{\l \in \Lambda} \left( M(\l) \times M(\l) \right) \rightarrow A
$$ is injective with image an $R$-basis of $A$.}
\item {(C2)}
{If $\l \in \Lambda$ and $S, T \in M(\l)$, we write $C(S, T) = C_{S, T}^{\l}
\in A$.  Then $*$ is an $R$-linear involutory anti-automorphism 
of $A$ such that
$(C_{S, T}^{\l})^* = C_{T, S}^{\l}$.}
\item {(C3)}
{If $\l \in \Lambda$ and $S, T \in M(\l)$ then for all $a \in A$ we have $$
a . C_{S, T}^{\l} \equiv \sum_{S' \in M(\l)} r_a (S', S) C_{S', T}^{\l}
\mod A(< \l),
$$ where  $r_a (S', S) \in R$ is independent of $T$ and $A(< \l)$ is the
$R$-submodule of $A$ generated by the set $$
\{ C_{S'', T''}^{\mu} : \mu < \l, S'' \in M(\mu), T'' \in M(\mu) \}
.$$}
\enddefinition

\remark{Remark 1.2.2}
We require the poset $\Lambda$ to be finite.  This is not part of
the original definition but problems occur in the general theory of
cellular algebras (see \cite{{\bf 13}, \S1.2}) 
and also in our later results if this hypothesis is omitted.
\endremark

There are many important examples of algebras in the mathematical
literature which turn out to be cellular, some of which we mention later.

\subhead 1.3 Tabular algebras \endsubhead

Table algebras and cellular algebras may be (usefully) amalgamated to form
``tabular algebras'' as follows.

\definition{Definition 1.3.1}
Let $\A = \zed[v, v^{-1}]$.  A {\it tabular algebra} is an
$\A$-algebra $A$, together with a table datum 
$(\Lambda, \Gamma, B, M, C, *)$ where:

\item{(A1)}
{$\Lambda$ is a finite poset.  For each $\l \in \Lambda$, 
$(\Gamma(\l), B(\l))$ is
a normalized table algebra over $\zed$ and
$M(\l)$ is a finite set (the set of ``tableaux'' of type $\l$).
The map $$
C : \coprod_{\l \in \Lambda} \left( M(\l) \times B(\l) \times M(\l)
\right) \rightarrow A
$$ is injective with image an $\A$-basis of $A$.  We assume
that $\Im(C)$ contains a set of mutually orthogonal idempotents 
$\{1_\e : \e \in {\Cal E}\}$ such that 
$A = \sum_{\e, \e' \in {\Cal E}} (1_\e A 1_{\e'})$ and such that for each
$X \in \Im(C)$, we have $X = 1_\e X 1_{\e'}$ for some $\e, \e' \in
{\Cal E}$.
A basis arising in
this way is called a {\it tabular basis}.  
}
\item{(A2)}
{If $\l \in \Lambda$, $S, T \in M(\l)$ and $b \in B(\l)$, we write
$C(S, b, T) = C_{S, T}^{b} \in A$.  
Then $*$ is an $\A$-linear involutory anti-automorphism 
of $A$ such that
$(C_{S, T}^{b})^* = C_{T, S}^{\overline{b}}$, where $\bar{\ }$ is the
table algebra anti-automorphism of $(\Gamma(\l), B(\l))$.
If $g \in \complex(v) \otimes_\zed \Gamma(\l)$ is such that 
$g = \sum_{b_i \in B(\l)} c_i b_i$ for some scalars $c_i$ 
(possibly involving $v$), we write
$C_{S, T}^g \in \complex(v)\otimes_\A A$ 
as shorthand for $\sum_{b_i \in B(\l)} c_i C_{S, T}^{b_i}$.  We write
$\bc_\l$ for the image under $C$ of $M(\l) \times B(\l) \times M(\l)$.}
\item{(A3)}
{If $\l \in \Lambda$, $g \in \Gamma(\l)$ and $S, T \in M(\l)$ then for all 
$a \in A$ we have $$
a . C_{S, T}^{g} \equiv \sum_{S' \in M(\l)} C_{S', T}^{r_a(S', S) g}
\mod A(< \l),
$$ where  $r_a (S', S) \in \Gamma(\l)[v, v^{-1}] = \A \otimes_\zed
\Gamma(\l)$ is independent of $T$ and of $g$ and $A(< \l)$ is the
$\A$-submodule of $A$ generated by the set $\bigcup_{\mu < \l} \bc_\mu$.}
\enddefinition

It is an easy consequence of these axioms that any table algebra (with
scalars extended to $\A$) is automatically a tabular algebra: set
$\Lambda$ and $M(\l)$ to be one-element sets, 
$(\Gamma, B)$ to be the table algebra in question, $*$ to be
the anti-automorphism of the table algebra and 
$C$ to be such that $C(m, b, m) = b$.  It is also clear that any
cellular algebra over $\A$ satisyfing the idempotent conditions in
(A1) is a tabular algebra: let $\Lambda, M, C,
*$ be as for cellular algebras, and let $(\Gamma(\l), B(\l))$ be the
one-dimensional table algebra spanned by the identity element.

\remark{Remark 1.3.2}
Note that if we apply $*$ to (A3), we obtain a condition (A$3'$) which
reads $$
C_{T, S}^{\overline{g}} . a^* \equiv \sum_{S' \in M(\l)} C_{T,
S'}^{\overline{g r_a (S', S)}}
\mod A(< \l)
.$$
\endremark

Next, we introduce an $\ba$-function (in the sense of Lusztig 
\cite{{\bf 22}, {\bf 23}, {\bf 24}, {\bf 25}})
associated to a tabular algebra $A$.

\definition{Definition 1.3.3}
Let $g_{X, Y, Z} \in \A$ be one of the structure constants for the tabular
basis $\Im(C)$ of $A$, namely $$
X Y = \sum_Z g_{X, Y, Z} Z
,$$ where $X, Y, Z \in \Im(C)$.   Define, for $Z \in \Im(C)$, $$
\afn(Z) = 
\max_{X, Y \in \Im(C)} \deg(g_{X, Y, Z})
,$$ where the degree of a Laurent polynomial is taken to be 
the highest power of $v$ occurring with nonzero coefficient.  We
define $\g_{X, Y, Z} \in \zed$ to be the coefficient of $v^{\afn(Z)}$ in $g_{X,
Y, Z}$; this will be zero if the bound is not achieved.
\enddefinition

Using the notion of $\ba$-function, we can now introduce ``tabular
algebras with trace''.

\definition{Definition 1.3.4}
A {\it tabular algebra with trace} is a tabular algebra in the sense
of Definition 1.3.1 which satisfies the conditions (A4) and (A5) below.

\item{(A4)}{Let $K = C_{S, T}^b$, $K' = C_{U, V}^{b'}$ and 
$K'' = C_{X, Y}^{b''}$ lie in $\Im(C)$.  Then the
maximum bound for $\deg(g_{K, K', K''})$
in Definition 1.3.3 is achieved if and only if $X = S$, $T = U$, $Y =
V$ and $b'' \in
\supp(bb')$ (see Definition 1.1.4).  If these conditions all hold and
furthermore $b = b' = b'' = 1$, we require $\g_{K, K', K''} = 1$.}
\item{(A5)}{There exists an $\A$-linear function $\t : A \ra \A$
(the {\it tabular trace}), such that $\t(x) = \t(x^*)$ for all $x \in
A$ and $\t(xy) = \t(yx)$ for all $x, y \in A$, that has the 
property that for every
$\l \in \Lambda$, $S, T \in M(\l)$, $b \in B(\l)$ and $X = C_{S,
T}^b$, we have $$
\t(v^{\afn(X)} X) = 
\cases 1 \mod v^{-1} \A^- & \text{ if } S = T \text{ and } b = 1,\cr
0 \mod v^{-1} \A^- & \text{ otherwise.} \cr
\endcases
$$  Here, $\A^- := \zed[v^{-1}]$.}
\enddefinition

We sketch a proof that the Hecke algebra of type $A$ is a tabular
algebra with trace.

\example{Example 1.3.5}
The Hecke algebra $\H = \H(A_{n-1})$ (over $\A$) 
of type $A_{n-1}$ is a tabular algebra.  The
table datum is an extension of the cell datum for $\H$ as a cellular
algebra which was given by Graham and Lehrer \cite{{\bf 8}, Example 1.2}.
In summary, the poset $\Lambda$ is the set of partitions of $n$, 
partially ordered by dominance.  The table algebras $(\Gamma(\l),
B(\l))$ are all trivial; that is, $B(\l) = \{1\}$ and $\Gamma(\l) =
\zed$.  The set $M(\l)$ is the set of standard tableaux of shape
$\l$.  Let $S, T$ be standard tableaux of shape $\l$.  The map $C$
takes the triple $(S, 1, T)$ to the Kazhdan--Lusztig basis element 
$C'_w$ via the Robinson--Schensted correspondence.  Since $1 \in \Im(C)$, the
idempotent condition is satisfied.  The map $*$ sends $C'_w$ to
$C'_{w^{-1}}$.

This algebra can be made into a tabular algebra with trace by taking
$\t(h)$ to be the coefficient of $T_1$ when $h$ is expressed as a
linear combination of the basis $\{v^{-\ell(w)} T_w: w \in
W(A_{n-1})\}$.  (The map $\tau$ was defined by Lusztig in \cite{{\bf 22},
\S1.4}.)  It follows from the definition of the $C'$-basis and the
properties of Kazhdan--Lusztig polynomials and Duflo involutions (see
the introduction to \cite{{\bf 23}}) that $\t(v^{\afn(C'_w)} C'_w)$ is $0 \mod
v^{-1} \A$ if and only if $w$ is not a Duflo involution.  Because all 
involutions are Duflo involutions in type $A$, axioms (A4) and (A5) follow
from \cite{{\bf 23}, Proposition 1.4} and properties of the
Robinson--Schensted map.
\endexample

Example 1.3.5 is trivial from the point of view of table algebras.
The next example, although much simpler in structure, is not.  We
introduce it since it will turn out that all tabular algebras of
finite rank with trace,
tensored over a suitable field, are isomorphic as abstract algebras to direct 
sums of algebras of the following form (Theorem 3.2.4).

\example{Example 1.3.6}
Let $n \in \enn$ and let $(\Gamma, B)$ be a normalized table algebra over 
$\zed$.  Let $M_n(\zed)$ be the ring of all $n \times n$ matrices with
integer coefficients.  Then the
algebra $\A \otimes_\zed M_n(\zed) \otimes_\zed \Gamma$ is tabular.  
For the table datum,
take $\Lambda$ to be a one-point set, $(\Gamma(\l), B(\l)) = (\Gamma, B)$,
$M(\l) = \{1, 2,
\ldots, n\}$, $C(a, b, c) = e_{ac} \otimes b$, where $e_{ac}$ is the
usual matrix unit, and $* : e_{ac} \otimes b \mapsto e_{ca} \otimes \bar{b}$.
Axioms (A1), (A2), (A3) are easily verified; in this case, the
mutually orthogonal idempotents in (A1) are the $n$ elements $e_{aa}
\otimes 1$.  The algebra can be made
into a tabular algebra with trace.  In this case, all basis elements
have $\ba$-value equal to $0$.  The product of two tabular basis
elements can be seen to be zero unless the conditions in axiom (A4)
hold, in which case the product of two basis elements is an integral
combination of other basis elements parametrized by elements of
$\supp(bb')$.  For (A5), we note that the map $\tau$ which takes
$e_{ac} \otimes b$ to $1$ if $a = c$ and $b = 1$, and to $0$ otherwise,
is a trace which satisfies the axiom.  (This works because the
operation of taking the coefficient of $1$ in a normalized table algebra is a
trace map; this is a consequence of axiom (T3).)
\endexample

The point of view taken in this paper is that the main objects of
interest are tabular bases, and not the abstract algebras they span.
This is in keeping with the philosophy behind table algebras.  A
finite dimensional normalized table algebra $(A, \BB)$ over $\real$ is
semisimple by Proposition 1.1.6, and therefore of little
interest as an abstract algebra $A$, although its basis $\BB$ has many
beautiful properties.  There may be some merit in considering
tabular algebras independently of their bases because, as K\"onig and
Xi have shown in \cite{{\bf 19}, {\bf 20}} and other papers, there 
is a basis-free approach to cellular algebras which leads to 
interesting results.

\head 2. Properties of tabular algebras \endhead

In this section, we investigate some of the consequences of the
tabular axioms.

\subhead 2.1 Relationship with cellular algebras \endsubhead

It turns out that many naturally occurring tabular algebras are also
known examples of cellular algebras.  The next result is a
sufficient criterion for a tabular algebra to be cellular.  This result
unifies several proofs already available in the literature.

\proclaim{Theorem 2.1.1}
Let $A$ be a tabular algebra of finite rank with table datum \newline
$(\Lambda, \Gamma, B, M, C, *)$; that is $|B(\l)| < \infty$ for each
$\l \in \Lambda$.  Suppose that, for some $R \geq \zed$ and for each
$\l \in \Lambda$, the algebra $R \otimes_\zed \Gamma(\l)$ is cellular 
over $R$ with cell datum $(\Lambda_\l, M_\l, C_\l, \bar{\ })$, where 
$\bar{\ }$ is the table algebra involution.  

Then $R \otimes_\zed A$ is cellular over $R \otimes_\zed \A$ with cell datum 
$(\Lambda', M', C', *)$, 
where $\Lambda' := \{(\l, \l') : \l \in \Lambda, \l' \in \Lambda_\l\}$ (ordered
lexicographically), $M'((\l, \l')) := M(\l) \times M_\l(\l')$ and $C((S,
s), (T, t))$ (where $(S, s), (T, t) \in M(\l) \times M_\l(\l')$) is
equal to $C_{S, T}^{C_\l(s, t)}$.
\endproclaim

\demo{Proof}
Axioms (C1) and (C2) follow immediately from the definitions and
axioms (A1) and (A2).  

To prove axiom (C3), let $\l \in \Lambda$ and 
let $C_\l(s, t)$ be a basis element of $\Gamma(\l)$ with $s, t \in
M_\l(\l')$.  Then by axiom (A3) we have, for any $a \in A$, $$
a . C_{S, T}^{C_\l(s, t)} \equiv \sum_{S' \in M(\l)} C_{S',
T}^{r_a(S', S) C_\l(s, t)} \mod A(< \l).
$$  Since $R \otimes_\zed \Gamma(\l)$ is cellular over $R$ with cell 
basis given by $C_\l$, it
follows by axiom (C3) applied to $R \otimes_\zed \Gamma(\l)$ that $$
r_a(S', S) C_\l(s, t) \equiv \sum_{s' \in M_\l(\l')} r'(S', S, s', s) 
C_\l(s', t)
\mod R \otimes \Gamma_\l (< \l'),
$$ where the $r'(S', S, s', s)$ are elements of $R \otimes_\zed \A$ which
are independent of $t$ (and, by axiom (A3), independent of $T$).
Axiom (C3) follows by tensoring over $R$.
\qed\enddemo

A good example of Theorem 2.1.1 concerns the Brauer algebra, which was shown
to be cellular by Graham and Lehrer \cite{{\bf 8}, \S4}.  (The reader
is referred to their paper for the definition.)  Here, we identify the 
parameter $\d$ with $v + v^{-1}$.

\example{Example 2.1.2}
Let $B(n)$ be the Brauer algebra (over $\A$) on $n$ strings.  Recall 
from \cite{{\bf 8},
\S4} that the algebra has an $\A$-basis consisting of certain triples $[S_1,
S_2, w]$ where $S_1, S_2$ are (arbitrary) involutions on $n$
letters with $t$ fixed points, and $w$ is an element of the symmetric
group $\Sy(t)$ on $t$ letters if $t > 0$, with $w = 1$ if $t = 0$.
The algebra has a table datum as follows.

Take $\Lambda$ to be the set of integers $i$ between $0$ and $n$ such that
$n - i$ is even, ordered in the natural way.  
If $\l = 0$, take $(\Gamma(\l), B(\l))$ to be the
trivial one-dimensional table algebra; otherwise, take $\Gamma(\l)$ to
be the group ring $\zed \Sy(t)$ with basis $B(\l) = \Sy(t)$ and 
involution $\overline{w}
= w^{-1}$.  Take $M(\l)$ to be the set of involutions on $n$ letters 
with $\l$ fixed points.  Take $C(S_1, w, S_2) = [S_1,
S_2, w]$; $\Im(C)$ contains the identity element.  
The anti-automorphism $*$ sends $[S_1, S_2, w]$ to $[S_2, S_1, w^{-1}]$.

Theorem 2.1.1 is applicable because $\Sy(t)$ is cellular over $\zed$;
this is immediate from setting $v = 1$ in Example 1.3.5.  We thus
recover Graham--Lehrer's cell datum for the Brauer algebra over $\A$.
\endexample

A technique similar to that used in Example 2.1.2 
may be applied to the case of the
partition algebra of \cite{{\bf 27}}; again the table algebras are
symmetric groups equipped with inversion as the involution.  This
recovers Xi's main result in \cite{{\bf 30}}.  
The generalized Temperley--Lieb algebra of type $H$, which is the
subject of \S5, is also covered by Theorem 2.1.1 (see \cite{{\bf 12},
\S3.3}); in this case we need $R$ to contain the roots of $x^2 - x - 1$.

\remark{Remark 2.1.3}
In fact, the algebra $B(n)$ can also be made into a tabular algebra with
trace.  The trace may be taken to have the property that $\tau(C_{S, S}^w)$ is 
$0$ if $w \ne 1$ and is $(v^{-1} + v^{-3})^k$ if $w = 1$ and 
$k = (n - \l)/2$, where $S \in M(\l)$.
\endremark

An interesting non-example of Theorem 2.1.1 which we shall mention
again later involves Jones' annular algebra, which we sketch below.  We
maintain the convention that $\d = v + v^{-1}$.  For the definition of 
the algebra, see \cite{{\bf 16}} and \cite{{\bf 8}, \S6}.

\example{Example 2.1.4}
Let $J_n$ be the Jones algebra (over $\A$) on $n$ strings.  Recall 
from \cite{{\bf 8},
\S6} that the algebra has an $\A$-basis consisting of certain triples $[S_1,
S_2, w]$ where $S_1, S_2$ are certain ``annular'' involutions on $n$
letters with $t$ fixed points, and $w$ is an element of the cyclic
group of order $t$ if $t > 0$, with $w = 1$ if $t = 0$.  The algebra 
has a table datum as follows.

Take $\Lambda$ to be the set of integers $i$ between $0$ and $n$ such that
$n - i$ is even, ordered in the natural way.  
If $\l = 0$, take $(\Gamma(\l), B(\l))$ to be the
trivial one-dimensional table algebra; otherwise, take $\Gamma(\l)$ to
be the group ring over $\zed$ of the cyclic group $\zed_\l$,
with basis $B(\l) = \zed_\l$ and involution $\overline{w}
= w^{-1}$.  Take $M(\l)$ to be the set of annular
involutions with $\l$ fixed points.  Take $C(S_1, w, S_2) = [S_1,
S_2, w]$, so that $\Im(C)$ contains the identity element.  
The anti-automorphism $*$ sends $[S_1, S_2, w]$ to
$[S_2, S_1, w^{-1}]$.
\endexample

Note that the table datum for $J_n$ is considerably simpler than the cell
datum for $J_n$ given in \cite{{\bf 8}, Theorem 6.15} and that it is
defined integrally in terms of a naturally occurring basis.  However,
we cannot apply Theorem 2.1.1 because group algebras of cyclic groups
are generally not cellular with respect to inversion.  Graham--Lehrer
thus need to use a more complicated involution than $*$ to establish
cellularity for $J_n$.  (See \cite{{\bf 20}, \S\S6--7} for more
details of the key role played by the involution in the structure of a
cellular algebra.)

\remark{Remark 2.1.5}
The algebra $J_n$ can also be made into a tabular algebra with trace.
The trace may be taken to have the property that $\tau(C_{S, S}^w)$ is 
$0$ if $w \ne 1$ and is $(v^{-1} + v^{-3})^k$ if $w = 1$ and 
$k = (n - \l)/2$, where $S \in M(\l)$.  The
calculation is similar to the calculations for generalized
Temperley--Lieb algebras which we perform in detail later.
We shall look at annular involutions more thoroughly in the context of
the affine Temperley--Lieb algebra in \S6.3.
\endremark

Two main advantages of considering a finite dimensional cellular
algebra as a tabular algebra are the following.  
First, tabular bases often arise out of natural constructions, such as bases of
Kazhdan--Lusztig type, whether or not one is motivated to define
tabular algebras.  (This also happens for cell bases, but to a lesser
extent.)  Secondly, there are many
examples of cellular algebras which do not have cell bases defined
over $\zed[v, v^{-1}]$.  Examples of these include Jones' annular
algebra above and the generalized Temperley--Lieb algebra of
type $H$ \cite{{\bf 12}}.  In each of these cases, the cell datum relies on
the fact that certain polynomials split over the ground ring, although
these polynomials do not split over $\zed$.  If these
algebras are considered as tabular algebras, they can be given natural
table data which are defined integrally.  

\subhead 2.2 Structure constants and the tabular trace \endsubhead

We next study some particular cases of axiom
(A3), analogous to the result \cite{{\bf 8}, Lemma 1.7}.

\definition{Definition 2.2.1}
Let $A$ be a tabular algebra with table datum $(\Lambda, \Gamma, B, M,
C, *)$.  Let $\l \in \Lambda$ and $S, T, U, V \in M(\l)$.  We define
$\lan T, U \ran \in \Gamma(\l)[v, v^{-1}]$ by the condition $$
C_{S, T}^1 C_{U, V}^{1} \equiv C_{S, V}^{\lan T, U \ran} \mod A(< \l)
.$$  If $b \in B(\l)$, we define $\lan T, U \ran_b \in \A$ to be the 
coefficient of $b$ in $\lan T, U \ran$.
\enddefinition

This is well-defined because of the following result.

\proclaim{Lemma 2.2.2}
Maintain the notation of Definition 2.2.1.  Then $$
C_{S, T}^b C_{U, V}^{b'}  \equiv C_{S, V}^{b \lan T, U \ran b'} \mod A(< \l)
$$ for any $b, b' \in B(\l)$, where $\lan T, U \ran$ is independent of
$b$, $b'$, $S$ and $V$.
\endproclaim

\demo{Proof}
First consider the expression $$
C_{S, T}^1 C_{U, V}^{1} \equiv C_{S, V}^{\lan T, U \ran} \mod A(< \l)
.$$  All basis elements $C_{S', V'}^b$ occurring on in the right hand
side must satisfy $V' = V$ by (A3) and $S' = S$ by (A$3'$).
Furthermore, $\lan T, U \ran$ is independent of $V$ by (A3) and
independent of $S$ by (A$3'$).  This shows that Definition 2.2.1 is
sound.

Starting from Definition 2.2.1, axiom (A3) implies that $$
C_{S, T}^1 C_{U, V}^{b'}  \equiv C_{S, V}^{\lan T, U \ran b'} \mod A(< \l)
$$ for any $b' \in B(\l)$ and then axiom (A$3'$) implies that $$
C_{S, T}^b C_{U, V}^{b'}  \equiv C_{S, V}^{b \lan T, U \ran b'} \mod A(< \l)
$$ for any $b \in B(\l)$, as required.
\qed\enddemo

\proclaim{Lemma 2.2.3}
Maintain the notation of Definition 2.2.1, with $b \in B(\l)$.
Denote the degree of the zero Laurent polynomial as $-\infty$ for notational
convenience, and suppose the tabular algebra $A$ satisfies axiom (A4).  
\item{\rm (i)}{We have $\deg \lan T, T
\ran_b \leq \deg \afn(C_{T, T}^1)$, with equality if and only if $b =
1$.  Moreover, $v^{\afn(C_{T, T}^1)}$ occurs in $\lan T, T
\ran_1$ with coefficient $1$.}
\item{\rm (ii)}{If $T \ne U$, we have $\deg \lan T, U
\ran_b < \deg \afn(C_{T, T}^1)$.}
\endproclaim

\demo{Proof}
We first prove (i).  Consider the expression $$
C_{T, T}^1 C_{T, T}^{\overline{b}} 
\equiv C_{T, T}^{\lan T, T \ran \overline{b}} \mod A(< \l)
.$$  By axiom (T3), the coefficient of $1$ in $\lan T, T \ran
\overline{b}$ is $\lan T, T \ran_b$, and hence the
coefficient of $C_{T, T}^1$ on the right hand side is also $\lan T, T
\ran_b$.  By axiom (A4),
$C_{T, T}^1$ appears on the right hand side with coefficient of
maximal degree if and only if $b = 1$, as $\supp(1\overline{b}) =
\{ \overline{b} \}$.  Furthermore, if $b = 1$, axiom (A4) guarantees
that the leading coefficient is $1$, as required.
Since the structure constants of $\Gamma(\l)$ do not involve $v$, (i)
follows.

To prove (ii), we consider the expression $$
C_{T, T}^1 C_{U, T}^{\overline{b}}
\equiv C_{T, T}^{\lan T, U \ran \overline{b}} \mod A(< \l)
.$$  Arguing as above, the coefficient of $C_{T, T}^1$ on the right
hand side is $\lan T, U \ran_b$.  On the other hand,
axiom (A4) shows that $C_{T, T}^1$ does not occur on the right hand side
with maximal degree, and (ii) follows.
\qed\enddemo

The tabular trace has the following key property.

\proclaim{Proposition 2.2.4}
Let $A$ be a (possibly infinite dimensional) tabular algebra with
trace.  Let $\l, \mu \in \Lambda$, $S, T \in M(\l)$, $U, V \in
M(\mu)$, $b \in B(\l)$ and $b' \in B(\mu)$.
Then $\t(C_{S, T}^b C_{U, V}^{b'}) = 0 \mod  v^{-1}\A^-$ unless $\l = \mu$,
$S = V$, $T = U$ and $\overline{b} = b'$.  If these conditions hold,
then $\t(C_{S, T}^b C_{U, V}^{b'}) = 1 \mod v^{-1}\A^-$.
\endproclaim

\demo{Proof}
Set $X := C_{S, T}^b$ and $X' := C_{U, V}^{b'}$.  Write $$
C_{S, T}^b C_{U, V}^{b'} = \sum_{X''} g_{X, X', X''} X''
,$$ where the sum is taken over tabular basis elements $X''$ so that
the $g_{X, X', X''}$ are structure constants with respect to the
tabular basis.  Now apply the tabular trace to both sides.  By axioms
(A4) and (A5), we see that $\t(g_{X, X', X''} X'')$ lies in $v^{-1}
\A^-$ if the $\afn$-function bound is not achieved.  Even if the bound
is achieved, meaning that $T = U$ and hence $\l = \mu$, 
axioms (A4) and (A5) imply that 
$\t(g_{X, X', X''} X'')$ will still lie in $v^{-1} \A^-$
unless $S = V$, $X'' = C_{S, S}^1$ and $1 \in \supp(bb')$.  
Axiom (T3) implies that the last condition happens if and only if 
$\overline{b} = b'$.  If all the conditions hold, we have $$
C_{S, T}^b C_{T, S}^{\overline b} \equiv C_{S, S}^{b \lan T, T \ran
\overline{b}} \mod A(< \l)
.$$  By Lemma 2.2.3 (i), the coefficient of $C_{S, S}^1$ in the right hand
side is a polynomial whose leading term is $v^{\afn(C_{T, T}^1)}$.
This occurs with coefficient $1$ because of the property of $\lan T, T
\ran_1$ described in Lemma 2.2.3 (i) and because the coefficient of $1 \in
\Gamma(\l)$ in $b \overline{b}$ is $1$.  Since the degree bound for
$g_{X, X', X''}$ is achieved for $X'' = C_{S, S}^1$, we have
$\afn(C_{T, T}^1) = \afn(C_{S, S}^1)$.
It follows that $\t(C_{S, T}^b C_{U, V}^{b'}) = 1 \mod v^{-1} \A$, as required.
\qed\enddemo

\proclaim{Theorem 2.2.5}
Let $A$ be a (possibly infinite dimensional) 
tabular algebra (over $\A$) with trace $\tau$ and table datum $(\Lambda,
\Gamma, B, M, C, *)$.  Then the map $(x, y) \mapsto \tau(xy^*)$
defines a symmetric, nondegenerate bilinear form on $A$ with the 
following properties.
\item{\rm (i)}{For all $x, y, z \in A$, $(x, yz) = (x z^*, y)$.}
\item{\rm (ii)}{The tabular basis is almost orthonormal with respect
to this bilinear form: whenever $X, X' \in \Im(C)$,
we have $$
(X, X') = 
\cases 1 \mod v^{-1} \A^- & \text{ if } X = X',\cr
0 \mod v^{-1} \A^- & \text{ otherwise.} \cr
\endcases$$}
\endproclaim

\demo{Proof}
Claim (i) follows from properties
of $*$, and claim (ii) is immediate from Proposition 2.2.4.  Nondegeneracy
follows easily from (ii).  Symmetry comes from axiom (A5):
$\t(xy^*) = \t((xy^*)^*) = \t(yx^*)$.
\qed\enddemo

\proclaim{Corollary 2.2.6}
A tabular algebra with trace is a symmetric algebra.
\endproclaim

\demo{Proof}
Maintain the notation of Theorem 2.2.5.  We need to show the existence
of a symmetric, associative and nondegenerate bilinear form.  The form
$f(x, y) := (x, y^*) = \t(xy)$ is clearly symmetric and associative,
and is nondegenerate because $(,)$ is nondegenerate.
\qed\enddemo

\subhead 2.3 Properties of the $\afn$-function \endsubhead

The $\afn$-function associated to a tabular algebra with trace has
similar properties to Lusztig's $\afn$-function from \cite{{\bf 22}}.  We
investigate some of these here.

\proclaim{Proposition 2.3.1}
Let $A$ be a tabular algebra with trace with table datum \newline
$(\Lambda, \Gamma, B, M,
C, *)$.  Let $\l \in \Lambda$, $b \in B(\l)$ and $S, U \in M(\l)$.  
The value of $\afn(C_{S, U}^b)$ depends only on $\lambda$, and not on 
$S \in M(\l)$, $U \in M(\l)$ or $b \in B(\l)$.
\endproclaim

\demo{Note}
We may use the notation $\afn(\l)$ in the sequel, with the
obvious meaning.
\enddemo

\demo{Proof}
By axiom (A4), we know that the coefficient of $C_{S, U}^b$ occurs
with maximal degree in the product $C_{S, T}^1 C_{T, U}^b$.  Now
consider the expression $$\eqalign{
C_{S, T}^1 C_{T, U}^{b} &\equiv C_{S, U}^{\lan T, T \ran b} \mod A(< \l)\cr
&\equiv \lan T, T \ran_1 C_{S, U}^b + \sum_{1 \ne b' \in B(\l)} C_{S,
U}^{\lan T, T \ran_{b'} b' b} \mod A(< \l). \cr
}$$  By Lemma 2.2.3 (i), the coefficient of $C_{S, U}^b$ occurring in the
first term has degree $\afn(C_{T, T}^1)$.  Again by Lemma 2.2.3 (i), the
coefficient of $C_{S, U}^b$ in each of the terms of the sum is either
zero or has degree strictly less than $\afn(C_{T, T}^1)$, because each
of the terms $\lan T, T \ran_{b'}$ is either zero or satisfies $\deg
\lan T, T \ran_{b'} < \afn(C_{T, T}^1)$.  
(Note that $v$ is not involved in the
expansion of $b'b$.)  Since $C_{S, U}^b$ occurs with maximal degree,
we must have $\afn(C_{S, U}^b) = \afn(C_{T, T}^1)$.  The claim follows.
\qed\enddemo

The next lemma is reminiscent of various results presented in \cite{{\bf 23}, \S1}.

\proclaim{Lemma 2.3.2}
Maintain the usual notation.  Let $A$ be a tabular algebra with trace 
and let $K_1, K_2, K_3$ be tabular basis elements.  Then $$
\g_{K_1, K_2, K_3^*} = \g_{K_2, K_3, K_1^*} = \g_{K_3, K_1, K_2^*}
,$$ where $\g$ is as in Definition 1.3.3.
\endproclaim

\demo{Proof}
In order for $\g_{K_1, K_2, K_3^*}$ 
to be nonzero, $K_3^*$ must appear with maximal
degree in the product $K_1K_2$.  If this happens, axiom (A4) requires 
(among other things)
that $K_1 = C_{S, T}^b$, $K_2 = C_{T, U}^{b'}$ and $K_3 = C_{U,
S}^{b''}$ for some $\l \in \Lambda$, $S, T, U \in M(\l)$ and $b, b',
b'' \in B(\l)$.  The same conditions on $S, T, U$ are necessary for
$\g_{K_2, K_3, K_1^*} \ne 0$ and for $\g_{K_3, K_1, K_2^*} \ne 0$, so
we can reduce consideration to the case where $S, T, U$ are as above.

Another requirement for $\g_{K_1, K_2, K_3^*} \ne 0$, again by axiom
(A4), is that $\overline{b''} \in \supp(bb')$.  If this condition is
met, then $\g_{K_1, K_2, K_3^*} \ne 0$; more precisely, $\g_{K_1, K_2,
K_3^*} = \k(\overline{b''}, bb')$ because $\lan T, T \ran_1$ has
leading coefficient $1$ by Lemma 2.2.3 (i).  Expansion of $bb' \in
\Gamma(\l)$ and axiom (T3) show that $\k(\overline{b''}, bb') = \k(1,
bb'b'')$.  Similar calculations show that $\g_{K_2, K_3, K_1^*} =
\k(1, b'b''b)$ and $\g_{K_3, K_1, K_2^*} = \k(1, b''bb')$.  Since
$\k(1, xy) = \k(1, yx)$ for all $x, y \in \Gamma(\l)$, we have $$
\k(1, bb'b'') = \k(1, b'b''b) = \k(1, b''bb')
,$$ and the claim follows.
\qed\enddemo

\proclaim{Corollary 2.3.3}
Let $A$ be a tabular algebra with trace 
and let $X, Y, Z \in \bc_\l$, as defined
in axiom (A2).  Then $$
\afn(Z) = 
\max_{X, Y \in \bc_\l} \deg(g_{X, Y, Z})
=
\max_{X, Y \in \bc_\l} \deg(g_{Y, Z, X})
=
\max_{X, Y \in \bc_\l} \deg(g_{Z, X, Y})
.$$
\endproclaim

\demo{Proof}
This follows from Proposition 2.3.1, Lemma 2.3.2
and the observation that as $X$ ranges over $\bc_\l$, so does $X^*$.
\qed\enddemo

\head 3. Asymptotic analogues of tabular algebras \endhead

Lusztig \cite{{\bf 26}} has developed a general method to send the parameter $v$ in
an $\A$-algebra to $\infty$ in a controlled way which leaves the
structure of the algebra essentially unchanged.  This approach is
valid if the algebra satisfies certain properties, which, as we shall
see in \S3.1, tabular algebras with trace satisfy.  This allows us to prove, in
\S3.2, a result giving the explicit structure of a tabular algebra
with trace over a suitable field.

\subhead 3.1 Cells \endsubhead

Lusztig's theory is designed to apply to
$\kyu(v)$-algebras equipped with bases with structure constants lying
in $\A$.  The basis is assumed to be compatible with a set of mutually
orthogonal idempotents as required by axiom (A1).

A slight problem to be overcome is that Lusztig's notion of cells
\cite{{\bf 26}, \S1.3}, which differs from the definition of
Kazhdan--Lusztig cells \cite{{\bf 18}, \S1}, is rather restrictive.  This
allows results to be proved for algebras without positivity of
structure constants, but is inconvenient here.  (All the examples of
tabular algebras in this paper have structure constants in $\enn[v,
v^{-1}]$, although we chose to omit this from the definition.)
Instead, we use the definition below; scrutiny of the results in
\cite{{\bf 26}, \S1} shows that they are all still valid, with the same
proofs, with this weaker notion of order.

\definition{Definition 3.1.1}
Let $A$ be a tabular algebra.  We introduce a relation, $\preceq$, on
the tabular basis by stipulating that $X' \preceq X$ if $X'$ appears
with nonzero coefficient in $KXK'$ for some tabular basis
elements $K, K'$.
\enddefinition

\remark{Remark 3.1.2}
If the structure constants of $A$ with respect to the tabular basis
lie in $\enn[v, v^{-1}]$, which is typical, the relation $\preceq$ is
automatically transitive.  We could also introduce one-sided (left or
right) versions of the relation $\preceq$ which would give rise below to left
cells and right cells as in \cite{{\bf 18}} or \cite{{\bf 5}, Definition 4.1};
we return to this briefly later.
\endremark

\proclaim{Proposition 3.1.3}
Let $A$ be a tabular algebra with table datum 
$(\Lambda, \Gamma, B, M, C, *)$ satisfying axiom (A4).  
Let $\preceq_t$ be the transitive
extension of the relation $\preceq$ of Definition 3.1.1.  The
relation $\sim$ on $\Im(C)$ defined by $Y \sim Z$ if and
only if $Y \preceq_t Z$ and $Z \preceq_t Y$ is an equivalence
relation.  The equivalence classes, known as $2$-cells, 
are parametrized by the elements of
$\Lambda$, where the class corresponding to $\l$ is $\bc_\l$.
\endproclaim

\demo{Proof}
The idempotent condition in axiom (A1) shows that $\sim$ is reflexive.

Let $Y = C_{T, U}^b$ and $Z = C_{V, W}^{b'} \in \bc_\l$; we will show
that $Y \sim Z$.  Now $$
YY^* = C_{T, U}^b C_{U, T}^{\overline{b}} \equiv C_{T, T}^{b \lan U, U
\ran \overline{b}} \mod A(< \l)
,$$ which, by Lemma 2.2.3 (i), contains $C_{T, T}^1$ with nonzero
coefficient of degree $\afn(\l)$.  There is a similar converse
statement: $C_{T, T}^1 C_{T, U}^b$ contains $Y$ with nonzero
coefficient.  This shows that $Y \sim C_{T, T}^1$.  Similarly, we have
$Z \sim C_{V, V}^1$.

Since $\sim$ is clearly symmetric and transitive, it is
an equivalence relation.  To finish the proof that $Y \sim Z$, we
observe that $C_{V, V}^1 \sim C_{S, S}^1$ for any $S, V \in M(\l)$.
This follows by consideration of the product $C_{S, V}^1 C_{V, V}^1
C_{V, S}^1$, which shows that $C_{S, S}^1 \preceq_t C_{V, V}^1$.

Because $\Lambda$ is partially ordered, axiom (A3) shows that the
equivalence classes of $\sim$ are no bigger than the sets $\bc_\l$
for fixed $\l$; these are therefore the equivalence classes.
\qed\enddemo

The definition of $\afn$-function in \cite{{\bf 26}} is in terms of
2-cells.  Let $L_\l$ be the $\A^-$-span 
of $\{ X: X \in \bc_\l\}$.  Then the $\afn$-function $\afn(Z)$, where
$Z \in \bc_\l$, is defined by Lusztig 
to be the smallest nonnegative integer $n$ such that $v^{-n} Z
L_\l \subseteq L_\l$, or $\infty$ if no such integer exists.  We also
define $A_\l$ to be the $\A$-submodule of $A$ spanned by $\bc_\l$.
This inherits an associative algebra structure from $A$ in the natural
way by setting $$
X X' = \sum_{X'' \in \bc_\l} g_{X, X', X''}X''
,$$ where the $g_{X, X', X''}$ are the structure constants for $A$.

\proclaim{Lemma 3.1.4}
For a tabular algebra with trace, 
Lusztig's definition of $\afn$-function agrees with Definition 1.3.3.
\endproclaim

\demo{Proof}
Let $X, Y, Z \in \Im(C)$ and let $\l \in \Lambda$ be such that $Z \in
\bc_\l$.  We have $$
\afn(Z) 
= \max_{X, Y \in \bc_\l} \deg(g_{X, Y, Z})
= \max_{X, Y \in \bc_\l} \deg(g_{Z, X, Y})
,$$ where the first equality is by axiom (A4) and the second is by
Corollary 2.3.3.  The claim follows.
\qed\enddemo

\subhead 3.2 Lusztig's properties $P_1$, $P_2$ and $P_3$ and
asymptotic tabular algebras \endsubhead

In order to send the parameter $v$ to $\infty$ in the correct way,
three properties ($P_1$, $P_2$ and $P_3$) are required of a quantum algebra.

\subsubhead Property $P_1$ \endsubsubhead

In \cite{{\bf 26}, \S1.4}, a basis $B$ is said to have property $P_1$ if (a) the
$\afn$-function takes finite values on $B$ and (b) for any 2-cell
$\bc_\l$ and any of the
orthogonal idempotents $1_\e$, the restriction of $\afn$ to $\bc_\l 1_\e$
is constant.

\proclaim{Lemma 3.2.1}
Let $A$ be a tabular algebra with trace and with table datum \newline
$(\Lambda, \Gamma, B, M, C, *)$.  Then $\Im(C)$ has property $P_1$.
\endproclaim

\demo{Proof}
The two notions of $\afn$-function agree by Lemma 3.1.4.  The
$\afn$-function is constant on 2-cells by Proposition 2.3.1, which
proves condition (b) of property $P_1$.  The $\afn$-function is
finite on any given 2-cell because $\lan T, T \ran \in \A$.
\qed\enddemo

Following \cite{{\bf 26}, \S1.4}, we write $\widehat{X} := v^{-\afn(X)} X$ for any
tabular basis element $X$.  The $\A^-$-submodule $A_\l^-$ of $A_\l$ is
defined to be generated by the elements $\{\widehat{X} : X \in
\bc_\l\}$.  We set $t_X$ to be the image of $\widehat{X}$ in $$
A_\l^\infty := {{A_\l^-} \over {v^{-1} A_\l^-}}
.$$  The latter is a $\zed$-algebra with basis $\{t_X : X \in
\bc_\l\}$ and structure constants $$
t_X t_{X'} = \sum_{X'' \in \bc_\l} \g_{X, X', X''} t_{X''}
,$$ where the $\g_{X, X', X''} \in \zed$ are as in Definition 1.3.3.
We also set $$
A^\infty := \bigoplus_{\l \in \Lambda} A_\l^\infty
;$$ this is a $\zed$-algebra with basis $\{t_X : X \in \Im(C)\}$.
It will turn out that, over a suitable field, $A^\infty$ is isomorphic to $A$.

\subsubhead Property $P_2$ \endsubsubhead

In \cite{{\bf 26}, \S1.5}, a basis $B$ with property $P_1$ is said to have property
$P_2$ if for any 2-cell $\bc_\l$, the $\zed$-algebra $A_\l^\infty$
admits a generalized unit.  This means that there is a subset $\D_\l$
of $\bc_\l$ that has the properties (a) that $t_D t_{D'} = \d_{D, D'}
t_D''$ whenever $D, D' \in \D_\l$ and (b) that for any $X \in \bc_\l$,
$t_X \in t_D A_\l^\infty t_{D'}$ for some (unique) $D, D' \in \D_\l$.

\proclaim{Lemma 3.2.2}
Let $A$ be a tabular algebra with trace and with table datum \newline
$(\Lambda, \Gamma, B, M, C, *)$.  Then $\Im(C)$ has property $P_2$.
\endproclaim

\demo{Proof}
Fix $\l$.  We define $\D_\l := \{C_{S, S}^1 : S \in M(\l)\}$.  If
$S \ne T \in M(\l)$, $D = C_{S, S}^1$ and $D' = C_{T, T}^1$, 
Lemma 2.2.3 (ii) shows that $t_D t_{D'} = 0$ in
$A_\l^\infty$.  The last condition of axiom (A4) shows that $t_D t_D =
t_D$.  This establishes part (a) of property $P_2$.  If $C_{S, T}^b
\in \bc_\l$, we set $D = C_{S, S}^1$ and $D' = C_{T, T}^1$ and then
Lemma 2.2.3 (i) gives the existence part (b) of property $P_2$;
uniqueness is by Lemma 2.2.3 (ii).
\qed\enddemo

\subsubhead Property $P_3$ \endsubsubhead

In \cite{{\bf 26}, \S1.6}, it is noted that there is a left $A$-module
structure on $A_\l$ given by $$
X . X' = \sum_{X'' \in \bc_\l} g_{X, X', X''} X''
,$$ where $X \in \Im(C)$ and $X', X'' \in \bc_\l$.  The same formula
with $X, X'' \in \bc_\l$ and $X' \in \Im(C)$ defines a right
$A$-module structure on $A_\l$.  One can then introduce a second
indeterminate, $v'$, giving rise to a $\zed[v', v^{\prime
-1}]$-algebra $A'$ which is the same
as $A$ except that $v$ is replaced by $v'$. We define a $\zed[v, v^{-1}, v', v^{\prime -1}]$-module
$\frak{A}_\l$ with basis $\{X: X \in \bc_\l\}$.  This is a left
$A$-module and a right $A'$-module, using the formulae above.
Following \cite{{\bf 26}, \S1.7}, we say a basis $B$ has property $P_3$ if
these two structures commute.

\proclaim{Lemma 3.2.3}
Let $A$ be a tabular algebra with trace and with table datum \newline
$(\Lambda, \Gamma, B, M, C, *)$.  Then $\Im(C)$ has property $P_3$.
\endproclaim

\demo{Proof}
Let $X$, $X'$ and $X''$ be a tabular basis element of
$A$, the basis element of $C_{U, V}^{b'}$ of $\frak{A}_\l$, 
and a tabular basis element of $A'$, respectively.
Then, using axioms (A3) and (A$3'$) and their notation, we obtain $$\eqalign{
(X . X') . ({X''}^*)
&= \sum_{U' \in M(\l)} C_{U', V}^{(r_X(U', U)(v)) b'} . ({X''}^*)\cr
&= \sum_{U', V' \in M(\l)} C_{U', V'}^{(r_X(U', U)(v)) b'
(\overline{r_{X''}(V', V)(v')})}\cr
&= \sum_{V' \in M(\l)} X . C_{U, V'}^{b' (\overline{r_{X''}(V', V)(v')})}\cr
&= X . (X' . ({X''}^*)).\cr
}$$  The superscripts in the sums are elements of $\Gamma(\l)[v,
v^{-1}, v', v^{\prime -1}]$.  The calculation works since the
structure constants of $\Gamma(\l)$ do not involve $v$ or $v'$.
\qed\enddemo

The fact that tabular algebras satisfy Lusztig's properties $P_1$,
$P_2$ and $P_3$ gives strong information about their structure.  In
\cite{{\bf 26}, \S1.8}, Lusztig defines a $\kyu(v)$-linear map 
$\Phi_\l : \kyu(v) \otimes_\A A \ra \kyu(v) \otimes_\zed A_\l^\infty$ which
satisfies $$
\Phi_\l(X) = \sum_{D \in \D_\l, Z \in \bc_\l} g_{X, D, Z} t_Z
$$ for $X \in \Im(C)$.  (There is no assumption that $X \in \bc_\l$.)
By \cite{{\bf 26}, Proposition 1.9 (b)}, this map is an algebra homomorphism
if the three properties are satisfied.  We can apply this result to
prove the following theorem.

\proclaim{Theorem 3.2.4}
Let $A$ be a tabular algebra of finite rank, with trace and with table
datum $(\Lambda, \Gamma, B, M, C, *)$.  
Let $k := \kyu(v)$.
\item{\rm(i)}{For any $\l \in \Lambda$, $A_\l \cong M_{|M(\l)|}(\zed)
\otimes_\zed \Gamma(\l)$ as $\zed$-algebras.}
\item{\rm(ii)}{There are $k$-algebra isomorphisms $$
k \otimes_\A A \cong k \otimes_\zed A^\infty \cong
\bigoplus_{\l \in \Lambda}
\left( k \otimes_\zed M_{|M(\l)|}(\zed) \otimes_\zed \Gamma(\l) \right).
$$}
\endproclaim

\demo{Proof}
We first prove (i) by showing that if $X = C_{S, T}^b$, the map
sending $t_X$ to $e_{S, T} \otimes b \in M_{|M(\l)|}(\zed) \otimes_\zed
\Gamma(\l)$ (where $e_{S, T}$ is a matrix unit) is a ring isomorphism.

Consider a product $t_X t_{X'}$ in $A_{\l}$.
Unless $X = C_{S, T}^b$ and $X' = C_{T, U}^{b'}$ for some $T$, the
degree bound in $XX'$ is not achieved and $t_X t_{X'}$ will be zero as
expected.  Otherwise, the tabular basis elements $X''$ which occur in
the product $XX'$ with maximal degree are, by axiom (A4), those of
form $C_{S, U}^{b''}$, where $b'' \in \supp(bb')$.  Since $$
C_{S, T}^b C_{T, U}^{b'} 
\equiv C_{S, U}^{b \lan T, T \ran b'} \mod A(< \l)
,$$ it follows from Lemma 2.3.2 and its proof that $\g_{X, X', X''} =
\k(b'', bb')$.  Similarly, the coefficient of $e_{S, U} \otimes b''$
in $(e_{S, T} \otimes b)(e_{T, U} \otimes b')$ is $\k(b'', bb')$.
This completes the proof of (i).

For the proof of (ii), we note that the isomorphism on the right in the
statement follows from (i), so we concentrate on proving $k \otimes_\A A
\cong k \otimes_\zed A^\infty$.  Let us define a homomorphism
$\Phi : \kyu(v) \otimes_\A A \ra \kyu(v) \otimes_\zed A^\infty$ via the direct
sum of the homomorphisms $\Phi_\l$ in each component, $A_\l^\infty$.

We claim that $\Phi$ is a monomorphism.  Let $X \in \Im(C)$, and consider
$\widehat{X}$.  Let $\l \in \Lambda$ (not necessarily such that $X \in
\bc_\l$).  
Then $$\Phi_\l(\widehat{X}) = \sum_{D \in \D_\l, Z \in \bc_\l} v^{-\afn(X)} 
g_{X, D, Z} t_Z.$$  By axiom (A4), Lemma 3.1.4 and the proof of Lemma
2.3.2, we see that all the terms on the right hand side
of this equation lie in $v^{-1} \A^{-} \otimes_\zed A_\l^\infty$ unless $X
\in \bc_\l$.  If $X \in \bc_\l$, there is exactly one term on the
right hand side for which this is not true, namely the one which
corresponds to the unique $D \in \D_\l$ with $t_X = t_X t_D$ (as in
Property $P_2$) and $Z = X$.  In this case, $\g_{X, D, Z} = 1$ 
and it follows that $$
\Phi_\l(\widehat{X}) - t_X \in v^{-1} \A^{-} \otimes_\zed A_\l^\infty
.$$  Considering all possibilities for $\l \in \Lambda$ gives 
$\Phi(\widehat{X}) - t_X \in v^{-1} \A^{-} \otimes_\zed A^\infty.$  
Suppose $x \in \ker(\Phi)$.  If $x \ne 0$, we may assume 
without loss of generality that 
the coefficients of $x$ with respect to the basis
$\{\widehat{X} : X \in \Im(C)\}$ lie in $\A^-$, but that not all of
them lie in $v^{-1} \A^-$.  The statement about $\Phi_\l(\widehat{X})
- t_X$ above shows that this cannot
happen and thus $\Phi$ is a monomorphism.  Since $A$ is of finite rank,
comparison of dimensions now completes the proof.
\qed\enddemo

\remark{Remark 3.2.5}
There is an interesting analogue of Theorem 3.2.4 for tabular algebras
of infinite rank, which involves the completion of an $\A^-$-form of
the algebra with respect to the $v^{-1}$-adic topology.  We omit the 
details for reasons of space.
\endremark

\head 4. Generalized Temperley--Lieb algebras of type $ADE$ \endhead

In the remaining sections of the paper, we look in detail at some
examples of tabular algebras with trace.  We illustrate the results
using generalized Temperley--Lieb algebras associated to Hecke
algebras of various kinds, starting in \S4 with the $ADE$ case.  The
structure of these algebras is well understood, and their combinatoric
properties have been analysed by Fan \cite{{\bf 5}} and others.  Apart from
the details of the tabular trace, most of the work required for the
verification of axioms (A1)--(A5) is done in the proofs of Fan's
results.

\subhead 4.1 Definitions \endsubhead

We start by defining the generalized Temperley--Lieb algebra $TL(X)$;
this coincides with the Temperley--Lieb algebra when $X$ is a Coxeter
graph of type $A_{n-1}$.

\definition{Definition 4.1.1}
Let $X$ be a Coxeter graph of type $A_n$, $D_n$ or $E_n$ for any $n
\in \enn$.  (We allow the long branch of a graph of type $E$ to be 
arbitrarily long.)  The associative, unital $\A$-algebra
$TL(X)$ is defined via generators $b_1, b_2, \ldots b_n$ corresponding
to the nodes of the graph and 
relations $$\eqalign{
b_i^2 &= [2] b_i, \cr
b_i b_j &= b_j b_i \quad \text{ if nodes $i$ and $j$ are not connected
in $X$,} \cr
b_i b_j b_i &= b_i \quad \text{ if nodes $i$ and $j$ are connected in $X$.}\cr
}$$  As usual, $[2] := v + v^{-1}$.
\enddefinition

Let $W(X)$ be the Coxeter group associated to $X$.
A product $w_1w_2\cdots w_n$ of elements $w_i\in W(X)$ is called
{\it reduced} if $$\ell(w_1w_2\cdots w_n)=\sum_i\ell(w_i).$$  We reserve
the terminology {\it reduced expression} for reduced products 
$w_1w_2\cdots w_n$ in which every $w_i \in S(X)$.

Call an element $w \in W(X)$ {\it complex} if it can be written 
as a reduced product $x_1 w_{ss'} x_2$, where $x_1, x_2 \in W(X)$ and
$w_{ss'}$ is the longest element of some rank 2 parabolic subgroup 
$\lan s, s'\ran$ such that $s$ and $s'$ do not commute.
Denote by $W_c(X)$ the set of all elements of $W(X)$
that are not complex.

For $w \in W_c(X)$, we define $b_w := b_{s_1} b_{s_2} \cdots b_{s_n}$
where $w = s_1 s_2 \cdots s_n$ is a reduced expression for $w$.  This
definition does not depend on the choice of reduced expression
\cite{{\bf 5}, \S2.2}.

\definition{Definition 4.1.2}
The set $\{b_w : w \in W_c(X)\}$ is an $\A$-basis for $TL(X)$.  We
call this the {\it monomial basis}.
\enddefinition

If $X$ is an arbitrary Coxeter graph, the generalized
Temperley--Lieb algebra may still be defined in a way which extends
Definition 4.1.1.  In this case, the algebra $TL(X)$ is the
quotient $\H(X)/J(X)$ of the usual $\A$-form of the Hecke algebra
$\H(X)$, where $J(X)$ is the two-sided ideal generated by the
Kazhdan--Lusztig basis elements $\{C'_w\}$ where $w$ is one of the
elements $w_{ss'}$ as above.  Such a generalized Temperley--Lieb
algebra is equipped with a canonical basis analogous to the
Kazhdan--Lusztig basis of the Hecke algebra.  The reader is referred
to \cite{{\bf 15}} for full details.

When $X$ is of type $ADE$, we may regard the monomial basis 
as natural in this context as it agrees with canonical basis for 
$TL(X)$ by \cite{{\bf 15}, Theorem 3.6}.  The purpose of \S4 is to show
that $TL(X)$ is a tabular algebra with trace having the monomial basis as its
tabular basis.

\subhead 4.2 Cellular structure and $\afn$-function in type $ADE$ \endsubhead

The following result is well-known and implicit in \cite{{\bf 5}}.

\proclaim{Proposition 4.2.1}
Let $X$ be a Coxeter graph of type $ADE$.  Then the algebra $TL(X)$ is
cellular with cell basis equal to the monomial basis and
anti-automorphism given by $* : b_w \mapsto b_{w^{-1}}$ for all $w \in W_c$.
\endproclaim

\demo{Proof}
The algebras in the statement have been shown to be of finite rank by 
Graham \cite{{\bf 7}, Theorem 7.1}.  
The decomposition of the monomial basis into cells (\idest 2-cells
in the sense of \S3.1) is described explicitly in \cite{{\bf 5}}; this provides
full details of the poset $\Lambda$ and the sets $M(\l)$.  In Fan's
terminology, the set $\Lambda$ is the set of two-sided cells ordered
by $\leq_{LR}$ (see \cite{{\bf 5}, Definition 4.1}), and $M(\l)$ can be
identified with the set of involutions in $W_c$ \cite{{\bf 5}, Theorem
4.4.4} which belong to the two-sided cell parametrized by $\l$.  This
works because a basis element is identifiable from the left cell and
the right cell which contain it, by \cite{{\bf 5}, Corollary 6.1.4}.  This
proves axiom (C1).  Let $C_{S, T}$ be the basis element in the same
left cell as $T$ and the same right cell as $S$.

Symmetry of the defining relations shows that the map
$*$ is an $\A$-linear anti-automorphism, and symmetry of the
definition of $W_c$ shows that it permutes the basis elements.
Symmetry of the definitions of left and right cells in \cite{{\bf 5},
Definition 4.1} shows that $*: C_{S, T} \mapsto C_{T, S}$.  This 
proves axiom (C2).

We can parametrize the basis
elements by $C_{S, T}$, where $S$ and $T$ are tableaux of the same
shape and correspond to ordered pairs of involutions in the same
two-sided cell by the above remarks.  The proof of \cite{{\bf 5}, Lemma
6.1.1} exhibits, for any $U \in M(\l)$, 
the existence of elements $a$ and $a'$ in $TL(X)$ (depending on $T$
and $U$ but not on $S$) which have the properties that $C_{S, T} a = C_{S, U}$ 
and $C_{S, U} a' = C_{S, T}$.  Axiom (C3) follows.
\qed\enddemo

Since the above algebras $TL(X)$ are cellular over $\A$ and the cell
basis contains the identity, they are trivially
tabular algebras.  The interesting aspect from our point of view is
that they are naturally tabular algebras with trace.

\proclaim{Lemma 4.2.2}
Let $x, y \in W_c(X)$.  Then $b_x b_y = [2]^m b_z$ for some $z \in
W_c$.  We have $m \leq \min(a(x), a(y)) \leq a(z)$, where $a$ is a certain
$\enn$-valued function depending only on the two-sided cell $\bc_\l$
containing $x$.  Furthermore, $m = \min(a(x), a(y))$ if and only if 
$b_x = C_{S, T}$ and $b_y = C_{T, U}$ for the same $T$.
\endproclaim

\demo{Proof}
We take $a$ to be the function defined in \cite{{\bf 5}, Definition
2.3.1}.  (It will turn out to coincide with our notion of
$\afn$-function, although this is not immediate.)  We have 
$\min(a(x), a(y)) \leq a(z)$ because $a(z) \geq a(x)$ and $a(z) \geq
a(y)$ by \cite{{\bf 5}, Corollary 4.2.2} and its dual.

The proof of \cite{{\bf 5}, Theorem 5.5.1} contains the proof of the other claims,
apart from the final statement which is an obvious consequence of
the argument given.
\qed\enddemo

\proclaim{Proposition 4.2.3}
Let $X$ be a Coxeter graph of type $ADE$.  Then the $\afn$-function
(in the sense of Definition 1.3.3) for the monomial basis of $TL(X)$
agrees with Fan's $a$-function in \cite{{\bf 5}} and it satisfies axiom (A4).
\endproclaim

\demo{Proof}
By \cite{{\bf 5}, Proposition 5.4.1}, we see that the structure constants
of $TL(X)$ are very simple: any two basis elements multiplied together
give a power of $[2]$ multiple of another basis element.  Consider the
equation $b_x b_y = [2]^m b_z$, 
and suppose that, for fixed $z$, $x$ and $y$ have been chosen to maximize the
number $m$ (meaning that $m = \afn(b_z)$).  
Lemma 4.2.2 now implies that $\afn(b_z) \leq a(z)$.  

We claim that $w$ and $w' \in W_c(X)$ can be chosen so that $b_w b_{w'} =
[2]^{a(z)} b_z$.  This will show that $a(z) \leq \afn(b_z)$.  Write
$b_z = C_{S, U}$ for some $S, U \in M(\l)$ and some $\l \in \Lambda$.
Choose $T$ such that $b_t := C_{T, T}$ is a product of commuting generators
$b_i$ in the two-sided cell corresponding to $\l$; this can be
arranged by \cite{{\bf 5}, Theorem 4.5.1}.
Define $w$ and $w'$ such that $b_w = C_{S, T}$ and $b_{w'} = C_{T,
U}$, so that $$
b_w b_{w'} = C_{S, T} C_{T, U} \equiv \lan T, T \ran C_{S, U} \mod A (< \l)
.$$  By \cite{{\bf 5}, Lemma 5.2.2}, $b_t b_t = [2]^{a(t)}b_t$, and thus
$\lan T, T \ran = [2]^{a(t)} \ne 0$.  Since the product of two basis
elements is a nonzero multiple of another, this means that $b_w b_{w'}
= [2]^{a(t)} b_z$, and that $w, w'$ and $z$ are in the same two-sided
cell.  Since the $a$-function is constant on two-sided cells, $a(t) =
a(z)$ and $a(z) \leq \afn(b_z)$.  This proves the first assertion.

Combining these results, we find that the only solutions of $b_x b_y =
[2]^m b_z$ for fixed $z$ and $m = \afn(z)$ require 
$m = \min(\afn(b_x), \afn(b_y)) = \afn(b_z) = a(z)$, and therefore
$b_x = C_{S, T}$, $b_y = C_{T, U}$ and $b_z = C_{S, U}$.  In this
case, $[2]^{\afn(z)}$ has leading coefficient $1$.
Axiom (A4) follows.
\qed\enddemo

\subhead 4.3 Cell modules and the tabular trace in types $ADE$ \endsubhead

Recall from \cite{{\bf 8}, Definition 2.1} that every cellular algebra
possesses a cell module for each $\l \in \Lambda$.  The parametrization of
the irreducible modules of a cellular algebra is in terms of the cell
modules.  Each cell module has a basis $\{C_S: S \in M(\l)\}$ and an
action of the algebra with the same structure constants as in axiom (C3).

\proclaim{Proposition 4.3.1}
Let $X$ be a Coxeter graph of type $ADE$.  Then the irreducible
modules for $TL(X)$ over the field $\kyu(\sqrt{[2]})$ are precisely
the cell modules.
\endproclaim

\demo{Proof}
By \cite{{\bf 5}, Theorem 5.6.1}, $TL(X)$ is semisimple over
$\kyu(\sqrt{[2]})$.  By \cite{{\bf 8}, Theorem 3.8}, this means that the
cell modules are the irreducible modules, and are pairwise nonisomorphic.
\qed\enddemo

\definition{Definition 4.3.2}
Let $X$ be a Coxeter graph of type $ADE$, and let $TL(X)$ be the
associated generalized Temperley--Lieb algebra, with cell datum
$(\Lambda, M, C, *)$, as in Proposition 4.2.1.  Let $\t'_\l$ be the
trace on the irreducible cell module corresponding to $\l \in
\Lambda$ via Proposition 4.3.1.  Then we define the trace $\t_\l$ on
$A$ via $\t_\l(x) := v^{-2\afn(\l)}\t'_\l(x)$.  (The $\afn$-function
is constant on $\bc_\l$ because it agrees with Fan's $a$-function by
Proposition 4.2.3, and Fan's $a$-function has this property.)
\enddefinition

\remark{Remark 4.3.3}
Notice that the cell module is defined with an $\A$-basis, so $\t'_\l(x)$ and
$\t_\l(x)$ both lie in $\A$ if $x$ is an $\A$-linear combination of
cell basis elements.  We will use this fact implicitly below without
further comment.
\endremark

\proclaim{Lemma 4.3.4}
Maintain the notation of Definition 4.3.2.  Let $Z = b_z$ be an element of
the cell basis, and let $\l \in \Lambda$ be arbitrary.  Then
$\t_\l(Z^*) = \t_\l(Z)$ and $$
\t_\l(v^{\afn(Z)}Z) = 
\cases 
1 \mod v^{-1} \A^- & \text{ if } Z = C_{S, S} \text{ for some } 
S \in M(\lambda),\cr
0 \mod v^{-1} \A^- & \text{ otherwise.} \cr
\endcases
$$
\endproclaim

\demo{Proof}
We use the notation of axiom (C3) as applied to cell modules.  For any
$T \in M(\l)$, we have $$
Z . C_{S, T}^{\l} \equiv \sum_{S' \in M(\l)} r_Z (S', S) C_{S', T}^{\l}
\mod A(< \l),
$$ which implies that $$
Z . C_S^{\l} = \sum_{S' \in M(\l)} r_Z (S', S) C_{S'}^{\l}
$$ in the cell module.  Applying the trace $\t_\l$ gives $$
\t_\l(v^{\afn(Z)} Z) = \sum_{T \in M(\l)} v^{-\afn(\l)} r_Z (T, T) 
.$$  Using axiom (A4), which holds by Proposition 4.2.3, we see that
$r_Z (T, T)$ is zero or has degree less than $\afn(\l)$ unless $Z =
C_{T, T}$.  In the latter case, $v^{-\afn(\l)}r_Z(T, T) - 1 \in v^{-1}
\A^-$.  The second claim follows.

To establish the first claim, we shall use the fact that the 
$\afn$-function is constant on two-sided cells.
Note also that if $a \in A(< \l)$, we have $\t_\l(a) = 0$.

Suppose $Z = C_{S, U}$, where $S, U \in M(\l)$.  Let $T \in M(\l)$.
As in the proof of Proposition 4.2.3, we have $\lan T, T \ran =
[2]^{\afn(\l)}$, so $$\eqalign{
\t_\l([2]^{\afn(\l)}Z) &= \t_\l(C_{S, T} C_{T, U})\cr
&= \t_\l(C_{T, U} C_{S, T})\cr
&= \lan U, S \ran \t_\l(C_{T, T})\cr
&= \lan S, U \ran \t_\l(C_{T, T})\cr
&= \t_\l([2]^{\afn(\l)}C_{U, S})\cr
&= \t_\l([2]^{\afn(\l)}Z^*).\cr
}$$  (Recall that $\lan , \ran$ is symmetric by \cite{{\bf 8},
Proposition 2.4 (i)}.)  It follows that $\t_\l(Z) = \t_\l(Z^*)$, as required.
\qed\enddemo

We can now prove the main result of \S4.

\proclaim{Theorem 4.3.5}
Let $X$ be a Coxeter graph of type $ADE$. Define a trace $\t : TL(X)
\ra \A$ by $\t(a) = \sum_{\l \in \Lambda} \t_\l(a)$.
Then $\t$ satisfies axiom (A5), and this makes 
$TL(X)$ with its canonical basis into a tabular algebra with trace.
\endproclaim

\demo{Proof}
It is clear from Lemma 4.3.4 that $$
\t(v^{\afn(Z)} Z) = 
\cases 1 \mod v^{-1} \A^- & \text{ if } Z = C_{S, S} \text{ for some }
S,\cr
0 \mod v^{-1} \A^- & \text{ otherwise.} \cr
\endcases
$$  The claim that
$\t(x) = \t(x^*)$ also follows from Lemma 4.3.4.  This establishes
axiom (A5) as the condition $b = 1$ is trivially true; the other axioms
follow from propositions 4.2.1 and 4.2.3.
\qed\enddemo

\head 5. Generalized Temperley--Lieb algebras of type $H$ \endhead

The generalized Temperley--Lieb algebras of type $H_n$ are an infinite
family of finite-dimensional algebras which arise as quotients of
(usually infinite dimensional) Hecke algebras associated to Coxeter
systems of type $H_n$.  The structure of this algebra was studied in
\cite{{\bf 12}} by means of a certain basis of diagrams with very
convenient properties.  It turns out that this basis is natural from
other perspectives: one of the main results of \cite{{\bf 14}} is that the
basis of diagrams agrees with the canonical basis for the algebra in
the sense of \cite{{\bf 15}}.  We show in \S5 that this basis has another
elegant property, namely that it is the tabular basis of a tabular algebra
with trace.  Unlike the situation in \S4, the table algebras involved
here are not trivial.

\subhead 5.1 The algebra of diagrams \endsubhead

We recall the definition of $TL(H_n)$, 
the generalized Temperley--Lieb algebra of type $H_n$, from \cite{{\bf 12}, \S1}.
This is based on a Coxeter system $X$ of type $H_n$ for $n \geq 2$ whose
Coxeter group $W(H_n)$ is given by generating involutions $\{s_i :
i \leq n\}$ and defining relations $$\eqalign{
s_i s_j &= s_j s_i \text{\quad if $|i - j| > 1$},\cr
s_i s_j s_i &= s_j s_i s_j \text{\quad if $|i - j| = 1$ and $\{i, j\}
\ne \{1, 2\}$},\cr
s_1 s_2 s_1 s_2 s_1 &= s_2 s_1 s_2 s_1 s_2.\cr
}$$  This is an infinite group for $n > 4$.  There is a corresponding
Hecke algebra with $\A$-basis $\{T_w: w \in W(H_n)\}$ and the usual
relations.

The $\A$-algebra $TL(H_n)$ is defined by the monomial basis
elements $b_i := b_{s_i}$ as follows:

\definition{Definition 5.1.1}
Let $n \in {\Bbb N} \geq 2$.  We define the associative, unital algebra
$TL(H_n)$ over $\A$ via generators $b_1, b_2, \ldots b_n$ and 
relations $$\eqalign{
b_i^2 &= [2] b_i, \cr
b_i b_j &= b_j b_i \text{\quad if \ $|i - j| > 1$},\cr
b_i b_j b_i &= b_i \text{\quad if \ $|i - j| = 1$ \ and \ $i, j > 1$},\cr
b_i b_j b_i b_j b_i &= 3 b_i b_j b_i - b_i
\text{\quad if \ $\{i, j\} = \{1, 2\}$}.\cr
}$$
\enddefinition

This can also be expressed in terms of the decorated tangles which
were defined in \cite{{\bf 11}}; for further elaboration and examples, the
reader is referred to \cite{{\bf 12}, \S2}.

A tangle is a portion of a knot diagram contained in a rectangle.  The
tangle is incident with the boundary of the rectangle only on the
north and south faces, where it intersects transversely.  The
intersections in the north (respectively, south) face are numbered
consecutively starting with node number $1$ at the western (\idest the
leftmost) end.
Two tangles are equal if there exists an isotopy of the plane carrying
one to the other such that the corresponding faces of the rectangle
are preserved setwise.
(We call the edges of the rectangular frame ``faces'' to avoid
confusion with the ``edges'' which are the arcs of the tangle.)

We extend the notion of a tangle so that each arc of the tangle may be
assigned a nonnegative
integer.  If an arc is assigned the value $r$, we represent this
pictorially by decorating the arc with $r$ blobs.

\definition{Definition 5.1.2}
A decorated tangle is a crossing-free tangle in which each arc is
assigned a nonnegative integer.  Any arc not exposed to the west face
of the rectangular frame must be assigned the integer $0$.

The category of decorated tangles, $\dt$, has as its objects the
natural numbers (not including zero).  The morphisms from $n$ to $m$
are the decorated tangles with $n$ nodes in the north face and $m$ in
the south.  The
source of a morphism is the number of points in the north face of the
bounding rectangle, and the target is the number of points in the
south face.  Composition of morphisms works by concatenation of the
tangles, matching the relevant south and north faces together.

Let $n$ be a positive integer.  The $\A$-algebra $\dt_n$ has as a
free $\A$-basis the morphisms from $n$ to $n$, where the multiplication
is given by the composition in $\dt$.

The edges in a tangle $T$ which connect nodes (\idest not the loops)
may be classified
into two kinds: propagating edges, which link a node in the north
face with a node in the south face, and non-propagating edges, which
link two nodes in the north face or two nodes in the south face.
\enddefinition

To explain how this relates to Coxeter systems of type $H$, we
recall from \cite{{\bf 12}, \S2.2} the notion of an $H$-admissible diagram.

\definition{Definition 5.1.3}
An $H$-admissible diagram with $n$ strands is an element of $\dt_n$
with no loops which satisfies the following conditions.

\item{\rm (i)}
{No edge may be decorated if all the edges are propagating.}
\item{\rm (ii)}
{If there are non-propagating edges in the diagram, then either there
is a decorated edge in the north face connecting nodes 1 and 2, or
there is a non-decorated edge in the north face connecting nodes $i$
and $i+1$ for $i > 1$.
A similar condition holds for the south face.}
\item{\rm(iii)}
{Each edge carries at most one decoration.}
\enddefinition

\topcaption{Figure 1} An $H$-admissible diagram for $n = 6$ \endcaption
\centerline{
\hbox to 2.638in{
\vbox to 0.888in{\vfill
        \includegraphics{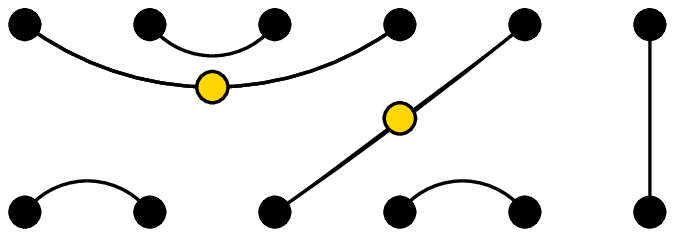}
}
\hfill}
}

\definition{Definition 5.1.4}
The algebra $\Delta_n$ (over a commutative ring with identity)
has as a basis the $H$-admissible
diagrams with $n$ strands and multiplication induced from that of
$\dt_n$ subject to the relations shown in Figure 2.
\enddefinition

\topcaption{Figure 2} Reduction rules for $\Delta_n$ \endcaption
\centerline{
\hbox to 1.041in{
\vbox to 1.916in{\vfill
        \includegraphics{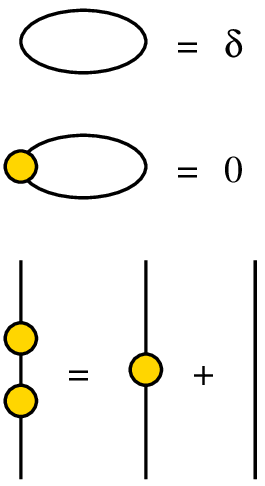}
}
\hfill}
}

By \cite{{\bf 12}, Lemma 2.2.4}, these rules give the free $\A$-module with
basis $\Delta_n$ the structure of an associative $\A$-algebra; we
identify $\d$ with $[2]$.  By \cite{{\bf 12}, Theorem 3.4.2}, we see that
$\Delta_n$ is isomorphic as an $\A$-algebra to $TL(H_{n-1})$.  The
isomorphism is an explicit one which identifies the generators $b_i$ with
certain $H$-admissible diagrams (see \cite{{\bf 12}, Proposition 3.1.2}).  
From now on, we will use this
identification implicitly and refer to both algebras as $TL(H_{n-1})$.

\subhead 5.2 Tabular structure in type $H$ \endsubhead

The main aim of \S5.2 is to show that the diagram basis for $TL(H_n)$
of \S5.1 is in fact a tabular basis, making the algebra into a tabular
algebra with trace.  This is particularly interesting
in light of \cite{{\bf 14}, Theorem 2.1.3}, which shows that the diagram
basis is also the canonical basis for $TL(H_n)$ in the sense of \cite{{\bf 15}}.

In order to understand the trace, it is convenient to make the
following definition.

\definition{Definition 5.2.1}
A trace diagram for $\dt_{n+1}$ consists of an isotopy class of 
$t$ (possibly zero) 
non-contractible, non-intersecting loops inscribed on a cylinder,
where the ends of the cylinder are labelled ``west'' and ``east''.
The westmost loop may optionally carry a single decoration, and the
integer $0 \leq t \leq n + 1$ must be such that $n + 1 - t$ is even.
We denote the set of trace diagrams for $\dt_{n+1}$ by $T(n+1)$, and
the free $\A$-module they span by $\A T(n+1)$.
\enddefinition

\definition{Definition 5.2.2}
Let $D$ be a decorated tangle with $n+1$ strands, subject to the
reduction rules in Figure 2.
Identify the north and south faces of $D$ to form a diagram inscribed on the 
surface of a cylinder.  The reduction rules in Figure 2 are clearly confluent
when applied to the cylinder, so the tangle can be reduced to an
element, $T(D)$, of $\A T(n+1)$.  An element in the support of $T(D)$
consists entirely of a disjoint union of $t$ (possibly zero) 
non-contractible loops on the cylinder, where the westmost loop may
optionally carry a single decoration.  (Note that $T(D)$ will be a linear 
combination of at most two trace diagrams, and will involve a power of 
$\d$ depending on the number of undecorated loops removed.  It is 
also possible that $T(D) = 0$.)  Furthermore, $n + 1 - t$ is
even, because there are $n+1$ strands in the diagram and the reduction
rules and isotopy preserve the parity of the number of intersections
with a west-east line.
\enddefinition

\proclaim{Lemma 5.2.3}
Let $f : T(n+1) \ra \A$ be any function; extend $f$ by linearity to
$\A T(n+1)$.  Then there is a 
trace $\t_f : TL(H_n) \ra \A$ defined by $\t_f(D) := f(T(D))$ and
extended linearly.
\endproclaim

\demo{Proof}
Let $D, D'$ be canonical basis elements.  It is clear from Definition
5.2.2 that $T(DD') = T(D'D)$; simply rotate the cylinder about its
axis by a half turn.  Since $f(T(DD')) = f(T(D'D))$, we see that
$\t_f(DD') = \t_f(D'D)$.  By linearity, $\t_f(xy) = \t_f(yx)$, as required.
\qed\enddemo

\definition{Definition 5.2.4}
Define $f : T(n+1) \ra \A$ as follows.  If $D \in T(n+1)$ carries a
decoration, set $f(D) = 0$.  
If $D$ is an element of $T(n+1)$ with $t$ non-contractible loops and no
decorations, set $f(D) = v^{t - (n+1)}$.
(Recall that $t - (n+1)$ is even and nonpositive.)  

Let $\t$ be the
trace on $TL(H_n)$ given by $\t_f$ for this value of $f$, as in Lemma 5.2.3.
\enddefinition

\proclaim{Theorem 5.2.5}
The algebra $TL(H_n)$ equipped with its canonical basis and 
the trace $\t$ of Definition 5.2.4 is a tabular algebra with trace.
\endproclaim

\demo{Proof}
We identify the canonical basis with a set of diagrams in $\dt_{n+1}$
in the usual way.  Let $\Lambda$ be the set of integers $t$ with $0
\leq t \leq n+1$ and $(n+1) - t$ even, ordered in the usual way.  

For $\l \in \Lambda$, let $(\G(\l), B(\l))$ be trivial if $\l = 0$ or $\l
= n+1$; for other values of $\l$, let $\G(\l) = \zed[x]/ \lan x^2 - x -
1 \ran$ and $B(\l) = \{1 + \lan x^2 - x - 1 \ran, x + \lan x^2 - x -
1 \ran \}$.  In all cases, the table algebra anti-automorphism is the
identity map.

Let $M(\l)$ be the set of possible configurations of non-propagating edges
in the north face of a canonical basis element with $\l$ propagating
edges.  The map $C$ produces
a basis element from the triple $(m, b, m')$ (where $m, m' \in M(\l),
b \in B(\l)$) as follows.  Turn the half-diagram corresponding to $m'$
upside down and place it below the half-diagram corresponding to $m$.
Join any free points in the bottom half to free points in the top half
so that they do not intersect.  If $b$ is the identity, leave all
propagating edges undecorated; if $b = x + \lan x^2 - x - 1\ran$,
decorate the westmost propagating edge.

The map $*$ is the linear extension of top-bottom inversion of the
diagrams.

Verification of axioms (A1) to (A3) is straightforward.  Notice that
the basis contains the identity element.

For axiom (A4), we claim that if $D \in \bc_\l$, we have 
$\afn(D) = a'(D) := ((n+1) - \l)/2$, \idest the
$\afn$-function evaluated at a diagram is half the number of
non-propagating edges in that diagram.  Let $D = C_{S, T}^b \in
\bc_\l$.  Direct computation shows that $C_{S, S}^1 D = [2]^{a'(D)} 
C_{S, T}^b$, so $\afn(D) \geq a'(D)$.  Conversely, the diagram
calculus shows that if $D'$ and $D''$ are canonical basis elements for
$TL(H_n)$, the number of loops formed in the product $D'D''$ is
bounded above both by $a'(D')$ and $a'(D'')$; this implies that the
structure constants appearing in $D'D''$ have degree bounded in the
same way.  Since $D$ can only appear
in a product $D'D''$ if $a'(D') \leq a'(D)$ and $a'(D'') \leq a'(D)$, we
have $\afn(D) \leq a'(D)$.  The claim follows.

The above argument also implies that the only way the
$\afn$-function bound can be achieved is if the three basis elements
$D', D'', D$ concerned come from the same $\bc_\l$.  The fact that
loops carrying a single decoration evaluate to zero means that the
bound can only be achieved if the pattern of edges at the bottom of
$D'$ is the same as the pattern of edges at the top of $D''$.  In this
case, we may set $D' = C_{S, T}^b$ and $D'' = C_{T, U}^{b'}$, and
properties of the diagram calculus give $D'D'' = [2]^{\afn(D)} C_{S,
U}^{bb'}$.  The assertions of axiom (A4) all follow easily.

Finally, we prove axiom (A5).  Consider a basis element $D$.  It is
clear by symmetry of the definitions that $\t(D) = \t(D^*)$, and thus
that $\t(x) = \t(x^*)$ for all $x \in A$.  To prove the other
requirements of the axiom,
we note that the diagram corresponding to $D$
has $2r$ non-propagating edges and $t$ propagating
edges, where $2r + t = n+1$ and $r = \afn(D)$.  To calculate $\t(D)$,
we inscribe $D$ on a cylinder as in Definition 5.2.2.  This forms a
number $r'$ of contractible loops, where we must have $r' \leq r$ since
each loop will include at least one non-propagating edge from each of
the north and the south faces of $D$.  In addition, each contractible
loop may include an even number of 
propagating edges from $D$ (e.g. $n = 5$ and $D
= b_2 b_3 b_4$).  Suppose the total number of propagating edges
involved in loops is $2t'$, where $2t' \leq t$.  We then have $$
\t(v^{\afn(D)}D) = \t(v^r D) = c v^{-r - 2t'} [2]^{r'} 
.$$ for some $c \in \enn$ (possibly zero).  In order for
$\t(v^{\afn(D)}D) \not\in v^{-1}\A^-$, we must have $r = r'$ and $t' =
0$.  For all these conditions to be met, we need 
$D = C_{S, S}^b$ for some $S \in M(r)$.  If $b \ne
1$, we have $\t(D) = 0$ by Definition 5.2.4.  If $b = 1$, we note that $$
\t(v^{\afn(D)}D) = (1 + v^{-2})^{\afn(D)} = 1 \mod v^{-1} \A^{-1}
.$$  Axiom (A5) follows, completing the proof.
\qed\enddemo
\head 6. The affine Temperley--Lieb algebra \endhead

The final examples of tabular algebras which we study here 
are the affine Tem\-per\-ley--Lieb
algebras, which turn out to be infinite dimensional tabular algebras.  
Affine Temperley--Lieb algebras 
are quotients of affine Hecke algebras which are
objects of considerable interest in representation theory.  Various
people have investigated the representation theory of the affine
Temperley--Lieb algebras, for example Martin and Saleur \cite{{\bf 28}},
the author and K. Erdmann \cite{{\bf 10}, {\bf 4}} and Graham and 
Lehrer \cite{{\bf 9}}.

It follows from the results of \cite{{\bf 6}, \S4} that the affine
Temperley--Lieb algebra has a faithful representation as an algebra of
diagrams.  By modifying this diagram basis slightly, we obtain a
tabular basis which is natural from the viewpoint of Kazhdan--Lusztig
theory.

\subhead 6.1 Diagrams for affine Temperley--Lieb algebras \endsubhead

We recall the graphical definition of the affine Temperley--Lieb
algebras from \cite{{\bf 4}, \S2.1}.  This can be given a more rigorous
treatment using categories analogous to our treatment of $TL(H_n)$ in
\S5.1, but this involves introducing technical definitions which are 
unnecessary for our purposes.  The reader is referred to \cite{{\bf 9},
\S1} for the categorical approach.

\definition{Definition 6.1.1}
An affine $n$-diagram, where $n \in {\Bbb Z}$ satisfies $n \geq 3$,
consists of two infinite horizontal rows of nodes lying at the points
$\{ {\Bbb Z} \times \{0, 1\}\}$ of ${\Bbb R} \times {\Bbb R}$,
together with certain curves, called edges, which satisfy 
the following conditions:

\item{\rm (i)}
{Every node is the endpoint of exactly one edge.}
\item{\rm (ii)}
{Any edge lies within the strip ${\Bbb R} \times [0, 1]$.}
\item{\rm (iii)}
{If an edge does not link two nodes then it is an infinite horizontal
line which does not meet any node.  Only finitely many edges are of
this type.}
\item{\rm (iv)}
{No two edges intersect each other.}
\item{\rm (v)}
{An affine $n$-diagram must be invariant under shifting to the left
or to the right by $n$.}
\enddefinition

By an isotopy between diagrams, we mean one which fixes the nodes and
for which the intermediate maps are also diagrams which are shift invariant.
We will identify any two diagrams which are isotopic to each other, so
that we are only interested in the equivalence classes of affine
$n$-diagrams up to isotopy.

Because of the condition {\rm (v)} in Definition 6.1.1,
one can also think of affine $n$-diagrams as
diagrams on the surface of a cylinder, or within an annulus, in a natural way.
We shall usually regard the diagrams as diagrams on the surface of
a cylinder with $n$ nodes on top and $n$ nodes on the bottom.
Under this construction, the top row of nodes
becomes a circle of $n$ nodes on one face of the cylinder, which we
will refer to as the top circle.  Similarly, the bottom circle of the
cylinder is the image of the bottom row of nodes.  We will use the
terms propagating and non-propagating to refer to edges as in
Definition 5.1.2.

\example{Example 6.1.2}
An example of an affine $n$-diagram for $n = 4$ is given in Figure 3.
The dotted lines denote the periodicity, and should be identified to
regard the diagram as inscribed on a cylinder.

\topcaption{Figure 3} An affine 4-diagram\endcaption
\centerline{
\hbox to 3.638in{
\vbox to 0.888in{\vfill
        \includegraphics{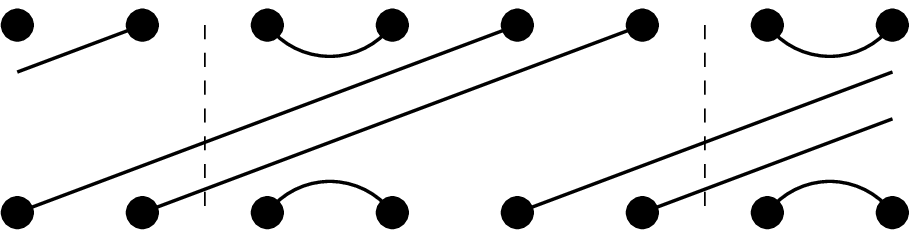}
}
\hfill}
}
\endexample

Two diagrams, $A$ and $B$ ``multiply'' in the following way, 
which was described in \cite{{\bf 6}, \S4.2}.
Put the cylinder for $A$ on top of the cylinder for $B$ and identify
all the points in the middle row.  This produces a certain (natural) number $x$
of loops.  Removal of these loops forms another diagram $C$ satisfying
the conditions in Definition 
6.1.1.  The product $AB$ is then defined to be $[2]^x C$.
It is clear that this defines an associative multiplication.

\definition{Definition 6.1.3}
The associative $\A$-algebra 
$\extln$ is the $\A$-linear span of all the affine 
$n$-diagrams, with multiplication given as above.
\enddefinition

\subhead 6.2 Generators and relations \endsubhead

We recall the presentation of $\extln$ by generators and relations.
For further explanation and examples, the reader is referred to
\cite{{\bf 4}, \S2.2}.

\definition{Definition 6.2.1}
Denote by $\bar{i}$ the congruence class of $i$ modulo $n$, taken from
the set $\text{\bf n} := \{1, 2, \ldots, n\}$.  We index the nodes in
the top and bottom circles of each cylinder by these congruence classes in the
obvious way.

The diagram $u$ of $\extln$ is the one satisfying the property that for
all $j \in \text{\bf n}$, the point $j$ in the bottom circle is connected
to point $\overline{j+1}$ in the top circle by a propagating edge taking
the shortest possible route.

The diagram $E_i$ (where $1 \leq i \leq n$) has
a horizontal edge of minimal length connecting $\bar{i}$ and 
$\overline{i+1}$ in each of the circles of the cylinder, and a
propagating edge connecting $\bar{j}$ in the top circle to $\bar{j}$ in
the bottom circle whenever $\bar{j} \ne \bar{i}, \overline{i+1}$.
\enddefinition

\proclaim{Proposition 6.2.2}
The algebra $\extln$ is generated by elements $$E_1, \ldots, E_n, u,
u^{-1}.$$  It is subject to the following defining 
relations: $$\eqalignno{
E_i^2 &= [2] E_i, & (1)\cr
E_i E_j &= E_j E_i, \quad \text{ if $\bar{i} \ne \overline{j \pm
1}$}, & (2) \cr
E_i E_{\overline{i\pm 1}} E_i &= E_i. & (3)\cr
u E_i u^{-1} &= E_{\overline{i+1}}, & (4)\cr
(u E_1)^{n-1} &= u^n . (u E_1). & (5)\cr
}$$
\endproclaim

\demo{Proof}
This is \cite{{\bf 10}, Proposition 2.3.7}.
\qed\enddemo

The algebra $\extll$ is closely related to the generalized
Temperley--Lieb algebra of type $\widehat{A}_l$, which is a quotient
of the Hecke algebra of type $\widehat{A}_l$ (see \cite{{\bf 9}, \S0}).

\proclaim{Proposition 6.2.3}
The algebra $\tln$ is the subalgebra of $\extln$ spanned by diagrams, $D$,
with the following additional properties: 

\item{\rm (i)}
{If $D$ has no horizontal edges, then $D$ is the identity diagram, in
which point $j$ in the top circle of the cylinder is connected to point
$j$ in the bottom circle for all $j$.}
\item{\rm (ii)}
{If $D$ has at least one horizontal edge, then the number of
intersections of $D$ with the line $x = i + 1/2$ for any integer $i$
is an even number.}

Equivalently, $\tln$ is the unital subalgebra of $\extln$ generated by
the elements $E_i$ and subject to relations (1)--(3) of Proposition 6.2.2.
\endproclaim

\demo{Proof}
See \cite{{\bf 10}, Definition 2.2.1, Proposition 2.2.3}.
\qed\enddemo

\subhead 6.3 Annular involutions \endsubhead

Let $D$ be an affine $n$-diagram associated to
the algebra $\extln$.  Until further notice,
we are only concerned with diagrams $D$ with $t > 0$ propagating
edges.  If $t > 0$, we can define the winding number $w(D)$ as follows.

\definition{Definition 6.3.1}
Let $D$ be as above.  
Let $w_1(D)$ be the number of pairs $(i, j) \in {\Bbb Z} \times {\Bbb Z}$
where $i > j$ and $\overline{j}$ in the bottom circle of
$D$ is joined to $\overline{i}$ in the top circle of $D$ by an edge
which crosses
the ``seam'' $x = 1/2$.  We then define $w_2(D)$ similarly but
with the condition that $i < j$, and we define $w(D) = w_1(D) - w_2(D)$.
\enddefinition

\definition{Definition 6.3.2 \cite{{\bf 8}, Lemma 6.2}}
An involution $S \in {\Cal S}_n$ is annular if and only if, for each
pair $i, j$ interchanged by $S$ $(i < j)$, we have
\item{\rm(i)}
{$S[i, j] = [i, j]$ and}
\item{\rm(ii)}
{$[i, j] \cap \text{\rm Fix} S = \emptyset$ or $\text{\rm Fix}S
\subseteq [i, j]$.}

We write $S \in I(t)$ if $S$ has $t > 0$ fixed points, and we write $S \in
\annn$ if $S$ is annular. In case $t=n$ we view the identity permutation as 
an annular involution.
\enddefinition

Using the concept of winding number, we have the following bijection
from \cite{{\bf 10}, Lemma 3.2.4}.

\proclaim{Proposition 6.3.3}
Let $D$ be a diagram for $\extln$ with $t$ propagating edges ($t > 0$).
Define $S_1, S_2 \in \annn \cap I(t)$ and 
$w \in {\Bbb Z}$ as follows.

The involution $S_1$
exchanges points $i$ and $j$ if and only if $i$ is connected to $j$ in
the top circle of $D$.  Similarly $S_2$ exchanges points $i$ and $j$ if
and only if $i$ is connected to $j$ in the bottom circle of $D$.
Set $w = w(D)$.

Then this procedure produces a bijection between diagrams $D$ with at
least one propagating edge and triples $[S_1, S_2, w]$ as above.
\endproclaim

\definition{Definition 6.3.4} Suppose $f(v) \in \zed[v, v^{-1}]$ is a 
Laurent polynomial, and
let $S_1, S_2 \in \annn \cap I(t)$ for $t > 0$. Then we write $f(S_1, S_2)
$ for the element of $\extln$ obtained from $f(v)$ by linear
substitution of $[S_1, S_2, k]$ for $v^k$. 
\enddefinition

Now consider the case where $t = 0$, as discussed in \cite{{\bf 4}, \S7.1}.

\definition{Definition 6.3.5}
Let $I(0)$ be the set of all permutations $S$ of $\Bbb Z$ which have the 
following properties:
\item{(a)}{for all $k \in {\Bbb Z}$, $S(n+k) = S(k) + n$;}
\item{(b)}{the image of $S$ in ${\Cal S}_n$ under reduction modulo $n$
is an annular involution with no fixed points.}
\enddefinition

With this assumption, there is
a one to one correspondence between $n$-diagrams with no propagating
edges and triples $[S_1, S_2, k]$ where $S_1, S_2 \in
I(0)$ and where $k \geq 0$. Here $S_1, S_2$ describe the top and the
bottom of the diagram as in Proposition 6.3.3, 
and $k$ is the number of infinite bands.  We
can extend Definition 6.3.4 to cover this situation as follows.

\definition{Definition 6.3.6} Suppose $f(x) \in \zed[x]$ is a 
polynomial, and let $S_1, S_2 \in I(0)$ (where $I(0)$ is as above).
Then we write $f(S_1, S_2)$ for the element of $\extln$ obtained from 
$f(x)$ by linear substitution of $[S_1, S_2, k]$ for $x^k$. 
\enddefinition

\subhead 6.4 Tabular structure \endsubhead

The basis of diagrams for $\extln$ can be made into a tabular basis
after some light modifications (which are trivial if $n$ is odd).  The
modifications to be made are determined by the properties of Chebyshev
polynomials.

\definition{Definition 6.4.1}
Let $\{U_k(x)\}_{k \in \enn}$ be the sequence of polynomials defined 
by the conditions $U_0(x) = 1$, $U_1(x) = x$ and the recurrence relation 
$U_{k+1}(x) = x U_k(x) - U_{k-1}(x)$.  
\enddefinition

Chebyshev polynomials are important in this context due to the
following result.

\proclaim{Proposition 6.4.2}
\item{\rm (i)}{The algebra $\zed[v, v^{-1}]$ with the involution
extending $\bar{v} = v^{-1}$ and basis $\{v^k : k \in \zed\}$ is a
normalized table algebra.}
\item{\rm (ii)}{The algebra $\zed[x]$ with the trivial involution
and basis $\{U_k : k \geq 0\}$ is a normalized table algebra.}
\endproclaim

\demo{Proof}
Part (i) is obvious, because $\zed[v, v^{-1}]$ is the group ring over
$\zed$ of the group of integers.

Part (ii) follows from the result $$
U_k(x) U_{k'}(x) = \sum_{i = 0}^{k} U_{k'-k + 2i}(x)
,$$ which is valid for $0 \leq k \leq k'$ and which can be established by an
easy induction.
\qed\enddemo

We are now ready to define the table datum for $\extln$.

\definition{Definition 6.4.3}
Let $n \geq 3$. 
Take $\Lambda$ to be the set of integers $i$ between $0$ and $n$ such that
$n - i$ is even; order $\Lambda$ in the usual way.  If $\l = 0$, take
$(\Gamma(\l), B(\l))$ to be the table algebra $\zed[x]$ of Proposition 6.4.2
(ii).  Otherwise, take $(\Gamma(\l), B(\l))$ to be the table algebra
$\zed[v, v^{-1}]$ of Proposition 6.4.2 (i).
For $\l \ne 0$, set $M(\l) := \annn \cap I(\l)$ as in Definition 6.3.2, and 
set $M(0) := I(0)$ as in Definition 6.3.5.
Take $C(S_1, f, S_2) = f(S_1, S_2)$ as in definitions 6.3.4 and 6.3.6;
note that $\Im(C)$ contains the identity element.  
The anti-automorphism $*$ corresponds to top-bottom reflection of the
diagrams.
\enddefinition

Note that if $n$ is odd, there is no cell $\l = 0$ and $\Im(C)$ above is the
same as the diagram basis.

We define the tabular trace for $\extln$ diagrammatically, following
similar lines to the argument in \S5.2.

\definition{Definition 6.4.4}
An affine trace diagram for $\extln$ consists of an isotopy class of 
$t$ (possibly zero) 
non-contractible, non-intersecting loops inscribed on the surface of a
torus.  We denote the set of affine trace diagrams for $\extln$ by 
$\wh{T}(n)$, and the free $\A$-module they span by $\A \wh{T}(n)$.
\enddefinition

\definition{Definition 6.4.5}
Let $D$ be an affine $n$-diagram inscribed on the surface of a cylinder.
Identify the top and bottom circles of $D$ to form a diagram inscribed on the 
surface of a torus.  
The seam $x = 1/2$ from Definition 6.3.2 maps to a non-contractible
circle on the torus which we shall also call the {\it seam}.
The resulting tangle can, by removal of contractible loops, be reduced to one
which consists entirely of a disjoint union of $t$ (possibly zero) 
non-contractible loops on the torus.  

There are two kinds of
non-contractible loops, which we call regular loops and exceptional
loops.  Exceptional loops become contractible in the solid torus and
arise (for example) from non-contractible loops in the cylinder on
which $D$ is inscribed.  Regular loops remain non-contractible in the
solid torus.  The non-intersection criterion guarantees that a tangle
arising from this construction cannot have both regular and
exceptional loops.
If there are no exceptional loops, the integer $n - t$ is even by an 
argument similar to that given for type $H$ in Definition 5.2.2.

Let $\wh{T}(D)$ be the element of $\A \wh{T}(n)$
arising from this construction.  
Note that $\wh{T}(D)$ will be a multiple of a single
affine trace diagram; the multiple is equal to $[2]^x$ where $x$ is 
the number of contractible loops removed.
\enddefinition

The proof of the following result is similar to that of Lemma 5.2.3,
but involves a torus rather than a cylinder.

\proclaim{Lemma 6.4.6}
Let $f : \wh{T}(n) \ra \A$ be any function; extend $f$ by linearity to
$\A \wh{T}(n)$.  Then there is a 
trace $\t_f : \extln \ra A$ defined by $\t_f(D) := f(\wh{T}(D))$ and
extended linearly.
\endproclaim

\definition{Definition 6.4.7}
Define $f : \wh{T}(n) \ra \A$ by its effect on $f(D)$ for $D \in
\wh{T}(n)$, as follows.  

Suppose first that $D$ has no exceptional
loops.  If all tangles isotopic
to $D \in \wh{T}(n)$ cross the seam of the torus, set $f(D) = 0$.  
If $D$ is isotopic to an element of $\wh{T}(n)$ with $t$
contractible loops such that none of the non-contractible loops crosses
the seam, set $f(D) = v^{t - n}$.
(Recall that in this case, $t - n$ is even and nonpositive.)  

Now suppose that $D$ has $k > 0$ exceptional loops.  Define $f(D) :=
v^{-n} \k(1, x^k)$, where $\k$ is the function in axiom (T3) applied
to the table algebra $\zed[x]$ of Proposition 6.4.2 (ii).

Let $\t$ be the
trace on $\extln$ given by $\t_f$ for this value of $f$, as in Lemma 6.4.6.
\enddefinition

\proclaim{Theorem 6.4.8}
The algebra $\extln$ equipped with the table datum of Definition 6.4.3
and the trace $\t$ of Definition 6.4.7 is a tabular algebra with trace.
\endproclaim

\demo{Proof}
The argument is reminiscent of the proof of Theorem 5.2.5, so we only
highlight differences.

Note that if $\l \ne 0$, $*$ sends $[S_1, S_2, k] \in \bc_\l$ to
$[S_2, S_1, -k]$.  However, if $\l = 0$, $*$ sends $[S_1, S_2, k]
\in \bc_\l$ to $[S_2, S_1, k]$.  Axiom (A2) now follows, and axioms
(A1) and (A3) are proved using routine arguments and the diagram calculus.

The $\afn$-value of a basis element may
be calculated by enumerating the half the number of non-propagating edges in
the diagram, not counting the non-contractible bands occurring in
elements of $\bc_{\l}$ for $\l = 0$.  Note that if $S \in M(0)$, we have $$
C_{S, S}^b C_{S, S}^{b'} = [2]^{n/2} C_{S, S}^{bb'}
,$$ where $n/2 = \afn(0)$.
The argument to prove axiom (A4) is now essentially the same as the one in
the proof of Theorem 5.2.5.  

The argument to establish axiom (A5) for a basis element $X \in
\Im(C)$ also follows easily from the
proof of Theorem 5.2.5, except in the case where $X \in \bc_0$.
If $X \not\in \bc_0$, we note that
if $X = C_{S, S}^b$, then when $X$ is inscribed on a torus as in
Definition 6.4.5, the non-propagating edges of $X$ must cross the seam
if and only if $b \ne 1$; the proof follows easily from this
observation.  If, on the other hand, $X = C_{S, S}^b$ and $X \in
\bc_0$, then $b = U_k(x)$ for some $k \in \zed^+$.  The function $f$
of Definition 6.4.7 then satisfies $f(X) = v^{-n} \k(1, U_k(x))$,
which evaluates to $v^{-n}$ if $k = 0$ and to $0$ otherwise.  The
verification of axiom (A5) is now routine, as $n = 2\afn(X)$.
\qed\enddemo

\head 7. Concluding remarks \endhead

The tabular basis in \S6 for $\extln$ is also connected to
Kazhdan--Lusztig theory.  If $A = \extln$ and $\BB$ is the tabular
basis for $A$, then there exists $\BB' \subset \BB$ and $A' \leq A$
such that $A'$ is the algebra $\tln$ of Proposition 6.2.3 and $\BB'$
is the canonical basis for $\tln$.  Furthermore, $\BB'$ turns out to be the
projection of the Kazhdan--Lusztig basis for $\H(\wh{A}_{n-1})$ to
$\tln$.  We call such a pair $(A', \BB') \leq (A, \BB)$ a {\it
sub-tabular algebra}.  There are many other natural examples of
sub-tabular algebras---for example, the Jones algebra
is a sub-tabular algebra of the Brauer algebra (compare examples 2.1.2
and 2.1.4)---and in a subsequent paper we shall consider how
generalized Temperley--Lieb algebras of types $B$ and $I$ also fit
into this framework.

An interesting question is whether any Hecke algebras of types other
than $A$ (recall Example 1.3.5), together with their Kazhdan--Lusztig bases,
are examples of tabular (or sub-tabular) algebras.  If so, it is
plausible that Lusztig's conjectures about Hecke algebras in \cite{{\bf 25},
\S10} could be interpreted in terms of sub-tabular algebras.

Finally, we note that 
Graham (unpublished) has investigated Markov traces for generalized 
Temperley--Lieb algebras which have some elegant properties.  It would be
interesting to know to what extent these are related to tabular traces.

\leftheadtext{}
\rightheadtext{}

\end